\documentclass[a4]{amsart}

\usepackage{amssymb}
\usepackage[all]{xy}
\usepackage{amsmath, amsfonts, amsthm , amssymb}

\input xypic

\newtheorem{theorem}{Theorem}[section]
\newtheorem{proposition}[theorem]{Proposition}
\newtheorem{lemma}[theorem]{Lemma}
\newtheorem{corollary}[theorem]{Corollary}

\theoremstyle{definition}
\newtheorem{definition}[theorem]{Definition}
\newtheorem{example}[theorem]{Example}

\newcommand{\lra}{\longrightarrow}

\newcommand{\Tor}{\ensuremath{\mathrm{Tor}}}

\newcommand{\Z}{\ensuremath{\mathcal{Z}}}

\newcommand{\F}{\ensuremath{\mathcal{F}}}
\newcommand{\C}{\ensuremath{\mathbb{C}}}
\newcommand{\CP}{\ensuremath{\mathbb{C}P^{\infty}}}
\newcommand{\uk}{\ensuremath{\mathrm{U}(K)}}
\renewcommand{\S}{\ensuremath{\mathcal{S}}}

\newcommand{\dr}[3]{\ensuremath{#1\stackrel{#2}
{\longrightarrow}#3}}
\newcommand{\ddr}[5]{\ensuremath{#1\stackrel{#2}
{\longrightarrow}#3\stackrel{#4}{\longrightarrow}#5}}

\newcounter{bean}
\newenvironment{letterlist}{\begin{list}{\rm ({\alph{bean}})}
      {\usecounter{bean}\setlength{\rightmargin}{\leftmargin}}}
      {\end{list}}

\newcommand{\namedright}[3]{\ensuremath{#1\stackrel{#2}
 {\longrightarrow}#3}}
\newcommand{\nameddright}[5]{\ensuremath{#1\stackrel{#2}
 {\longrightarrow}#3\stackrel{#4}{\longrightarrow}#5}}
\newcommand{\namedddright}[7]{\ensuremath{#1\stackrel{#2}
 {\longrightarrow}#3\stackrel{#4}{\longrightarrow}#5
  \stackrel{#6}{\longrightarrow}#7}}

\newcommand{\larrow}{\relbar\!\!\relbar\!\!\rightarrow}
\newcommand{\llarrow}{\relbar\!\!\relbar\!\!\larrow}
\newcommand{\lllarrow}{\relbar\!\!\relbar\!\!\llarrow}

\newcommand{\lnamedright}[3]{\ensuremath{#1\stackrel{#2}
 {\larrow}#3}}

\newcommand{\llnamedright}[3]{\ensuremath{#1\stackrel{#2}
 {\llarrow}#3}}

\newcommand{\llnamedddright}[7]{\ensuremath{#1\stackrel{#2}
 {\llarrow}#3\stackrel{#4}{\llarrow}#5
  \stackrel{#6}{\llarrow}#7}}

\newcommand{\lllnameddright}[5]{\ensuremath{#1\stackrel{#2}
 {\lllarrow}#3\stackrel{#4}{\lllarrow}#5}}

\newcommand{\qqed}{\hfill\Box}

\begin{document}
\title[The homotopy type of the complement of an arrangements]{The
homotopy type of the complement of a coordinate subspace
arrangement}
\author{Jelena Grbi\'{c} and Stephen Theriault}
\address{Department of Mathematical Sciences,
         University of Aberdeen, Aberdeen AB24 3UE, United Kingdom}
\email{jelena@maths.abdn.ac.uk}
\email{s.theriault@maths.abdn.ac.uk}

\date{}

\begin{abstract}
The homotopy type of the complement of a complex coordinate subspace
arrangement is studied by fathoming out the connection between
its topological and combinatorial structures. A family of
arrangements for which the complement is homotopy equivalent to a
wedge of spheres is described. One consequence is an application in
commutative algebra: certain local rings are proved to be Golod,
that is, all Massey products in their homology vanish.
\end{abstract}

\maketitle \tableofcontents

\newpage

\section{Introduction}

In this paper we study connections between the topology of the
complements of certain complex arrangements, and  algebraic and
combinatorial objects associated to them.

Let
\begin{equation*}
\ensuremath{\mathcal{A}}=\{L_{1},\ldots ,L_{r}\}
\end{equation*}
be a \textit{complex subspace arrangement} in
$\ensuremath{\mathbb{C}}^{n}$,
that is, a finite set of complex linear subspaces in $\ensuremath{\mathbb{C}}%
^{n}$. For such an arrangement $\ensuremath{\mathcal{A}}$, define
its
\textit{support} $|\ensuremath{\mathcal{A}}|$ as $|\ensuremath{\mathcal{A}}%
|=\bigcup_{i=1}^{r}L_{i}\subset \ensuremath{\mathbb{C}}^{n}$ and its \textit{%
complement} $U(\ensuremath{\mathcal{A}})$ as
\begin{equation*}
U(\ensuremath{\mathcal{A}})=\ensuremath{\mathbb{C}}^{n}\backslash |%
\ensuremath{\mathcal{A}}|.
\end{equation*}

Arrangements and their complements play a pivotal role in many
constructions of combinatorics, algebraic and symplectic geometry,
etc.; they also arise as configuration spaces for different
classical mechanical systems. Special problems connected with
arrangements and their complements arise in different areas of
mathematics and mathematical physics. The multidisciplinary nature
of the subject results in ongoing theoretical improvements,
a constant source of new applications and
the penetration of new ideas and techniques in each of the component
research areas. It is the
interplay of methods from seemingly disparate areas that makes the
theory of subspace arrangements a vivid and appealing field of
research.

In the study of arrangements it is important to get a detailed
description of the topology of their complements, including
properties such as homology groups, cohomology rings, homotopy type,
and so on. In this paper we are concerned with the homotopy type of
the complement of a complex coordinate subspace arrangement. A {\it
complex coordinate subspace} of $\mathbb{C}^{n}$ is given by
\begin{equation*}
L_{\sigma }=\{(z_{1},\ldots ,z_{n})\in \mathbb{C}^{n}\,|\
z_{i_{1}}=\cdots =z_{i_{k}}=0\}
\end{equation*}
where $\sigma =\{i_{1},\ldots ,i_{k}\}$ is a subset of $[n]=\{1,\ldots ,n\}$%
, allowing us to define a {\it complex coordinate subspace arrangement} $\mathcal{C%
}\ensuremath{\mathcal{A}}$ in $\ensuremath{\mathbb{C}}^{n}$ as a
family of coordinate subspaces $L_{\sigma }$ for $\sigma \subset
\lbrack n]$. The main topological space we study, naturally
associated to the complex coordinate subspace arrangement
$\mathcal{C}\ensuremath{\mathcal{A}}$, is
the complement $U(\mathcal{C}\ensuremath{\mathcal{A}})$ in $%
\ensuremath{\mathbb{C}}^{n}$. Our results are obtained by studying
the topological and combinatorial structures of
$U(\mathcal{C}\mathcal{A})$ with the help of commutative and
homological algebra, combinatorics and homotopy theory.

It has been known for some time that hyperplane arrangements have a
torsion free cohomology ring. Recently it was proved~\cite{S}
that after suspending the complement of a hyperplane arrangement it
becomes homotopy equivalent to a wedge of spheres. The case of
complex coordinate subspace arrangements is much more complicated.
Already at the cohomology level, there is a more intricate
structure. The Buchstaber-Panov formula for
$H^*(U(\mathcal{C}\ensuremath{\mathcal{A}}))$~\cite{BP}
detects torsion in special cases, implying that even stably $U(\mathcal{C}%
\ensuremath{\mathcal{A}})$ cannot always be homotopy equivalent to a
wedge of spheres. That makes the question of when the complement of
a coordinate subspace arrangement is homotopy equivalent to a wedge
of spheres more difficult and therefore more interesting. The main
goal of this paper is to describe a family of coordinate subspace
arrangements for which the complement is homotopy equivalent to a
wedge of spheres.

The basic connections between the topology, combinatorics and
commutative algebra of coordinate subspace arrangements are
established as follows.

Let $K$ be a simplicial complex on the vertex set $[n]$. We shall
consider only complexes that are finite, abstract simplicial
complexes represented by
their collection of faces. Every simplicial complex $K$ on the vertex set $%
[n]$ defines a complex arrangement of coordinate subspaces in $%
\ensuremath{\mathbb{C}}^{n}$ via the correspondence
\begin{equation*}
K\ni \sigma \mapsto \mathrm{span}\{e_{i}:i\not\in \sigma \}
\end{equation*}
where $\left\{ e_{i}\right\} _{i=1}^{n}$ is the standard basis for $%
\ensuremath{\mathbb{C}}^{n}$. Equivalently, for each simplicial
complex $K$ on the set $[n]$, we associate the complex coordinate
subspace arrangement
\begin{equation*}
\mathcal{C}\ensuremath{\mathcal{A}}(K)=\{L_{\sigma }|\,\sigma
\not\in K\}
\end{equation*}
and its complement
\begin{equation}
\label{uk} U(K)=\mathbb{C}^{n}\backslash \bigcup_{\sigma \not\in
K}L_{\sigma }.
\end{equation}
On the other hand, to $K$ and a commutative ring $R$ with unit there
is an associated algebraic object, the {\it Stanley-Reisner ring}
$R[K]$, also known in the literature as the {\it face ring} of $K$.
Denote by $R[v_{1},\ldots ,v_{n}]$ the graded polynomial algebra on
$n$ variables where $\deg (v_{i})=2$ for each $i$ over $R$. The
Stanley-Reisner ring of a simplicial complex $K$ on the vertex set
$[n]$ is the quotient ring
\begin{equation*}
R[K]=R[v_{1},\ldots ,v_{n}]/\mathcal{I}_{K}
\end{equation*}
where $\mathcal{I}_{K}$ is the homogeneous ideal generated by all
square free monomials $v^{\sigma }=v_{i_{1}}\cdots v_{i_{s}}$ such
that $\sigma =\{v_{i_{1}},\ldots v_{i_{s}}\}\not\in K$.

Coming back to topology and following the Buchstaber-Panov
approach~\cite{BP} to toric topology, there are another two topological
spaces associated to a simplicial complex $K$ and its
Stanley-Reisner ring $R[K]$. The first space arises as a topological
realisation of the Stanley-Reisner ring. It is the
\textit{Davis-Januszkiewicz space} $DJ (K),$ whose cohomology ring
is isomorphic to the Stanley-Reisner ring $R[K]$. The
Davis-Januszkiewicz space maps by an inclusion into the classifying
space of the $n$-dimensional torus. The homotopy fibre of this inclusion
can be
identified with another torus space, the \textit{moment-angle complex} $%
\ensuremath{\mathcal{Z}}_{K},$ which has as a deformation retract the
complement $U(K)$ of the complex coordinate subspace
arrangement~\cite[8.0]{BP}.
Different models of $DJ (K)$ and $\ensuremath{\mathcal{Z}}%
_{K}$ as well as their additional properties will be addressed later
on in Section~\ref{sec:mainobjects}. These homotopic identifications
show that the problem of determining the homotopy type of the
complement of complex coordinate subspace arrangements is equivalent
to determining the homotopy type of the moment-angle complex
$\ensuremath{\mathcal{Z}}_{K}$. To do this we need to closely
examine the homotopy fibration sequence
\begin{equation*}
\ensuremath{\Z_K\stackrel{}
{\longrightarrow}DJ(K)\stackrel{\mathrm{incl}}{\longrightarrow}BT^n}.
\end{equation*}
The main technique employed for understanding of this filtration is
Mather's Cube Lemma~\cite{M}, which relates homotopy pullbacks
and homotopy pushouts in a cubical diagram. This is applied
iteratively as $K$ is built up one face at a time, in a prescribed
order. An analysis of the component homotopy fibration and
cofibration sequences produces our main result,
Theorem~\ref{shifted_wedge}(see below).

To find a suitable simplicial complex $K$ whose $U(K)$ will be
homotopy equivalent to a wedge of spheres, we first look at its
cohomology ring. As $U(K)$ is homotopy equivalent to $\Z_K$, this is
the same as looking at the cohomology ring of $\Z_K$. The integral
cohomology of $\ensuremath{\mathcal{Z}}_K$ has been calculated in
\cite[7.6 and 7.7]{BP}. If $\mathcal{Z}_{K}$ is to be homotopy
equivalent to a wedge of spheres then we need to consider simplicial
complexes $K$ for which all Massey products in
$H^*(\ensuremath{\mathcal{Z}}_K)$ vanish. That will not imply that
$\ensuremath{\mathcal{Z}}_K$ is itself homotopic to a wedge of
spheres but at least on the cohomological level there will be no
obstructions to that claim. Combinatorists, from their point of
view, have studied simplicial complexes and associated to them
certain Tor algebras that correspond to the cohomology of
$\ensuremath{\mathcal{Z}}_K$ as in our case. They have determined
several classes of complexes for which it can be shown that all
Massey products in associated Tor algebras vanish.

One such class is of shifted complexes. A simplicial complex $K$ is
\textit{shifted} if there is an ordering on the vertex set such that
whenever $\sigma$ is a simplex of $K$ and $v^{\prime}<v$, then
$(\sigma - v)\cup v^{\prime}$ is a simplex of $K$. Gasharov, Peeva
and Welker~\cite{GPW} showed that when $K$ is a shifted complex,
then all Massey products in $H^*(\ensuremath{\mathcal{Z}}_K)$ are
trivial. In this
case we obtain much stronger result by determining the homotopy type of $%
\ensuremath{\mathcal{Z}}_K$.

\begin{theorem}
\label{shifted_wedge} Let $K$ be a shifted complex. Then $U(K)$ is
homotopy equivalent to a wedge of spheres.
\end{theorem}

Previously, the only known cases of simplicial complexes $K$ for
which the complement $U(K)$ has the homotopy type of a wedge of
spheres occurred when $K$ was a disjoint union of $n$ vertices. When
$n=2$ or $n=3$, these are classical results of low dimensional
topology, while the general case was proved by the
authors~\cite{GT}. The result in Theorem~\ref{shifted_wedge} is much
more general. Notice for example that any full $k$-dimensional
skeleton of the standard simplicial complex on $n$ vertices $\Delta^{n}$
is a shifted complex.

Let $\ensuremath{\mathcal{F}}_t$ be the family of simplicial
complexes $K$ for which the moment-angle complex
$\ensuremath{\mathcal{Z}}_K$ has the property that $\Sigma^t
\ensuremath{\mathcal{Z}}_K$ is homotopy equivalent to a wedge of
spheres. Our next theorem describes the influence that
combinatorial operations on simplicial complexes have with respect to $%
\ensuremath{\mathcal{F}}_t$.

\begin{theorem}
\label{families} Let $K_{1}\in \ensuremath{\mathcal{F}}_{t}$ and $K_{2}\in %
\ensuremath{\mathcal{F}}_{s}$ for some non-negative integers $t$ and~$s$.
The effect on family membership of the simplicial
complex $K$ resulting from the following operations on $K_1$ and $K_2$ is:

\begin{enumerate}
\item  the disjoint union of simplicial complexes:\newline
if $K=K_{1}\coprod K_{2}$, then $K\in \ensuremath{\mathcal{F}}_{m}$
where $m=\max \{t,s\} $;
\item  gluing along a common face:\newline
if $K=K_{1}\bigcup_{\sigma }K_{2}$, then $K\in
\ensuremath{\mathcal{F}}_{m}$ where $\sigma $ is a common face of
$K_{1}$ and $K_{2}$ and ${m=\max \{t,s\}}$;

\item the join of simplicial complexes:\newline
if $K=K_{1}\ast K_{2}$, then $K\in \ensuremath{\mathcal{F}}_{m}$
where $m=\max \{t,s\}+1$.
\end{enumerate}
\end{theorem}

As a corollary we specify the operations on simplicial complexes for
which $\ensuremath{\mathcal{F}}_{0}$ is closed.

\begin{corollary}
\label{familiescor}
Let $K_{1}$ and $K_{2}$ be simplicial complexes in $\ensuremath{%
\mathcal{F}}_{0}$. Then $\ensuremath{\mathcal{F}}_{0}$ is closed for
the following operations on simplicial complexes:

\begin{enumerate}
\item  the disjoint union of simplicial complexes,\newline
$K=K_{1}\coprod K_{2}\in \ensuremath{\mathcal{F}}_{0}$;

\item  gluing along a common face,\newline
$K=K_{1}\bigcup_{\sigma }K_{2}\in \ensuremath{\mathcal{F}}_{0}$, where $%
\sigma $ is a common face of $K_{1}$ and $K_{2}$.
\end{enumerate}
\end{corollary}

The information we have obtained on complex subspace arrangements
has an application in commutative algebra. Let $R$ be a local ring.
One of the fundamental aims of commutative algebra is to describe
the homology ring of $R$, that is $\Tor _R(k,k)$, where $k$ is a
ground field. The first step in understanding
$\ensuremath{\mathrm{Tor}}_R(k,k)$ is to obtain information about
its Poincar\'{e} series $P(R)$, more specifically, whether  $P(R)$
is a rational function. A far reaching contribution to this problem
was made by Golod.
A local ring $R$ is \textit{Golod} if all Massey products in $%
\ensuremath{\mathrm{Tor}}_{k[v_1,\ldots,v_n]}(R,k)$ vanish.
Golod~\cite{Go} proved that if a local ring is Golod, then its
Poincar\'{e} series represents a rational function and it is
determined by $P(\Tor_{k[v_1,\ldots,v_n]}(R,k))$. Although being
Golod is an important property, not many Golod rings are known.
Using our results on the homotopy type of the complement of
a coordinate subspace arrangement, we are able to use homotopy theory
to gain some insight into these difficult homological-algebraic
questions. The main results are as follows.

\begin{theorem}
\label{ineq} For a simplicial complex $K$,
\begin{equation*}
P(k[K])\leq \frac{t(1+t)^n}{t-P(H^*(U(K)))}.
\end{equation*}
Equality is obtained when $k[K]$ is Golod.
\end{theorem}

\begin{theorem}
\label{eq} If $K\in \mathcal{F}_{0}$, then $k[K]$ is a Golod ring.
\end{theorem}

Combining Theorems~\ref{ineq} and~\ref{eq}, we obtain the following
result.

\begin{corollary}
For a simplicial complex $K\in \mathcal{F}_{0}$,
\begin{equation*}
P(k[K])= \frac{t(1+t)^n}{t-P(H^*(U(K)))}.
\end{equation*}
\end{corollary}

To close, let us remark that all the techniques used in this paper
can be also applied to real and quaternionic coordinate subspace
arrangements by changing the ground ring from complex numbers to
real, quaternion numbers respectively. In those cases
Theorem~\ref{shifted_wedge} describes the homotopy type of the
complement of real, quaternionic coordinate subspace arrangements.
For real arrangements instead of torus spaces and
$\mathbb{C}P^{\infty}$, we look at spaces with an action of
$\mathbb{Z}/2$ (also considered as $S^0$) and
$\mathbb{R}P^{\infty}$, respectively; while in the case of
quaternionic arrangements we deal with $S^3$ spaces and
$\mathbb{H}P^{\infty}$.

The disposition of the paper is as follows.
Section~\ref{sec:mainobjects} catalogues the main objects of study
and states various properties they satisfy.
Sections~\ref{sec:prelim} through~\ref{sec:ZKtype} build up to and
deal with the primary focus of the paper,
Theorem~\ref{shifted_wedge}. Sections~\ref{sec:prelim}
through~\ref{sec:fatwedge} establish the preliminary homotopy
theory. Included are identifications of the homotopy types of
various pushouts, a review of homotopy actions, the general
statement of Mather's Cube Lemma and a finer analysis of a special
case involving homotopy actions, and several properties of the fat
wedge. Section~\ref{sec:fibreg} considers a particular pattern of
successive inclusions of one coordinate subspace into another which
we term a regular sequence. Such a sequence need not always exist,
but when it does we show there is a measure of control over the
homotopy types of the successive homotopy fibres obtained from
including the coordinate subspaces into the full coordinate space
$X_{1}\times\cdots\times X_{n}$. Section~\ref{sec:regseqs} gives
conditions guaranteeing the existence of regular sequences, which
are based on the properties of a shifted complex.
Section~\ref{sec:ZKtype} puts together all the material in
Sections~\ref{sec:prelim} through~\ref{sec:regseqs} to prove
Theorem~\ref{shifted_wedge}. At this point, the class of simplicial
complexes for which $\mathcal{Z}_{K}$ is homotopy equivalent to a
wedge of spheres includes the shifted complexes.
Section~\ref{sec:topext} shows that there are other simplicial
complexes $K$ which have $\mathcal{Z}_{K}$ homotopy equivalent (or
stably homotopy equivalent) to a wedge of spheres by proving
Theorem~\ref{families} and Corollary~\ref{familiescor}. Finally,
Section~\ref{sec:algebra} turns to commutative algebra considering
Golods rings and their properties, and proves Theorems~\ref{ineq}
and~\ref{eq}.

\noindent\textbf{Acknowledgements.} The authors would like to thank
Professors Victor Buchstaber and Taras Panov for their stimulating
work,  as well as for their helpful suggestions and kind
encouragement. The first author would also like to thank Professor
Volkmar Welker for explaining to her the connection between
combinatorics and arrangements and for making it possible for her to
visit the University of Marburg for a week.

\section{The main objects: their definitions and properties}
\label{sec:mainobjects}

As mentioned in the introduction the main objective of this paper is
the study of arrangements and their complements from topological
point of view. To pass from the combinatorial concept of
arrangements to a topological one, we use different topological
models associated to simplicial complexes $K$ and their algebraic
counterparts, the Stanley-Reisner rings $\mathbb{Z}[K]$ (or the face
rings) of $K$.

The purpose of this section is to present the main objects which we are
going to use and to set the notation. We rely heavily on
constructions in toric topology introduced and studied  by
Buchstaber and Panov~\cite{BP}.

\subsection{The Davis-Januszkiewicz space} The topological realisation of
the Stanley-Reisner ring $\mathbb{Z}$ is called the Davis-Januszkiewicz
space $DJ(K)$. The first model of $DJ(K)$ is a Borel-type construction
due to Davis and Januszkiewicz~\cite{DJ}. For our
purposes we use another model of $DJ(K)$ given by
Buchstaber-Panov~\cite{BP}. In what follows, we identify the
classifying space of the circle $S^1$ with the infinite-dimensional
projective space $\CP$, and therefore the classifying space $BT^n$
of the $n$-torus with the $n$-fold product of $\CP$. For an
arbitrary subset $\sigma\subset[n]$, define the $\sigma$-power of
$BT$ as
\[
BT^{\sigma}=\{(x_1,\ldots, x_n)\in BT^n\, |\ x_i=* \text{ if }
i\notin\sigma\}.
\]
\begin{definition}
Let $K$ be a simplicial complex on the index set $[n]$. The {\it
Davis-Januszkiewicz} space is given as the cellular subcomplex
\[
DJ(K)=\bigcup_{\sigma\in K}BT^\sigma\subset BT^n.
\]
\end{definition}

Buchstaber and Panov justified the name of this topological model by
proving the following.
\begin{proposition}[Buchstaber-Panov~\cite{BP}]
The cohomology of $DJ(K)$ is isomorphic to the Stanley-Reisner ring
$\mathbb{Z}[K]$. Moreover, the inclusion of cellular complexes
$i\colon DJ(K)\lra BT^n$ induces the quotient epimorphism
\[
i^*\colon\mathbb{Z}[v_1,\ldots,v_n]\lra \mathbb{Z}[K]
\]
in cohomology.
\end{proposition}
Recently, Notbohm-Ray~\cite{NR} showed that the Davis-Januszkiewicz
spaces are uniquely determined, up to homotopy equivalence, by
their cohomology ring. This implies that all models of
Davis-Januszkiewicz spaces are mutually homotopy equivalent.

\subsection{The moment-angle complex}
Realise the torus $T^n$ as a subspace of $\C^n$
\[
T^n=\big\{(z_1,\ldots,z_n)\in \mathbb{C}^n \,|\ |z_i|=1, \text{ for
} i=1,\ldots ,n\, \big\}
\]
contained in the unit polydisc
\[
(D^2)^n=\big\{(z_1,\ldots,z_n)\in \mathbb{C}^n \,|\ |z_i|\leq 1,
\text{ for } i=1,\ldots ,n\, \big\}.
\]
For an arbitrary subset $\sigma\subset[n]$, define
\[
B_\sigma=\big\{(z_1,\ldots,z_n)\in (D^2)^n \,|\ |z_i|=1 \quad
i\notin \sigma\big\}.
\]
\begin{definition}
Let $K$ be a simplicial complex on the index set $[n]$. Define the
{\it moment-angle complex} $\Z_K$ by
\[
\Z_K=\bigcup_{\sigma\in K}B_\sigma \subset (D^2)^n.
\]
\end{definition}
Observe that since each $B_\sigma$ is invariant under the action of
$T^n$, the moment-angle complex $\Z_k$ is a $T^n$-space. Buchstaber
and Panov showed that the moment-angle complex is another
topological model of the Stanley-Reisner ring $\mathbb{Z}[K]$ by
proving that the $T^n$-equivariant cohomology $H^*_{T^n}(\Z_K)$ is
isomorphic to $\mathbb{Z}[K]$.

The following description of the moment-angle complex $\Z_K$
together with its relation to the complement of an arrangement plays
the pivotal role in our approach to determine the homotopy type
of the complement of a complex coordinate subspace arrangement.

\begin{proposition}[Buchstaber-Panov~\cite{BP}]
The moment-angle complex $\Z_K$ is the homotopy fibre of the
embedding
\[
i\colon DJ(K)\lra BT^n.
\]
\end{proposition}

Recall from \eqref{uk} that $\uk$ denotes the complement of the
complex coordinate subspace arrangement associated to a simplicial
complex $K$.
\begin{theorem}[Buchstaber-Panov~\cite{BP}]
\label{retract} There is an equivariant deformation retraction
\[
\dr{\uk}{}{\Z_K}.
\]
\end{theorem}

\subsection{The cohomology of moment-angle complexes and shifted complexes}
Theorem~\ref{retract} insures that the homotopy type of the
complement $\uk$ of a complex coordinate subspace arrangement can be
obtained by finding the homotopy type of the moment-angle complex
$\Z_K$.

In our study of the homotopy type of $\uk$ we specialised by asking
for which simplicial complexes $K$ the complement $\uk$ is homotopy
equivalent to a wedge of spheres. To begin we look at the cohomology
ring of $\Z_K$ finding those simplicial complexes for which there is
no cohomological obstruction for $\Z_K$ to be homotopic to a wedge
of spheres. Buchstaber and Panov~\cite{BP} described the cohomology
algebra of $Z_K$ by proving that there is an isomorphism
\[
 H^*(\Z_K;k)\cong
\Tor_{k[v_1,\ldots,v_n]}(k[K],k).
\]
as graded algebras.

\begin{definition}
\label{Goloddef}
The Stanley-Reisner ring $k[K]$ is {\it Golod} if all Massey
products in $\Tor_{k[v_1,\ldots,v_n]}(k[K],k)$ vanish.
\end{definition}
This definition provides a class of rings $k[K]$ for which $Z_K$
might be homotopic to a wedge of spheres. Although being Golod is an
interesting property of a ring, there are not many examples of Golod
rings. The one that is going to be of use for us comes from
combinatorics.
\begin{definition}
A simplicial complex $K$ is {\it shifted} if there is an ordering on
its set of vertices such that whenever $\sigma\in K$ and $v'<v$,
then $(\sigma - v)\cup v'\in K$.
\end{definition}
Notice that any full $i$-th skeleton $\Delta^i(n-1)$ of the standard
simplicial complex $\Delta^{n-1}$ (also denoted by $\Delta(n)$) on
$n$ vertices is shifted.
\begin{proposition}[Gasharov, Peeva and Welker~\cite{GPW}]
If $K$ is shifted, then its face ring $k[K]$ is Golod.
\end{proposition}

\section{Preliminary homotopy decompositions}
\label{sec:prelim}

The purpose of this section is to identify the homotopy
type of several pushouts. We begin by stating Mather's
Cube Lemma~\cite{M}, which relates homotopy pullbacks and
homotopy pushouts in a cubical diagram.

\begin{lemma}
   \label{cube}
   Suppose there is a homotopy commutative diagram
   
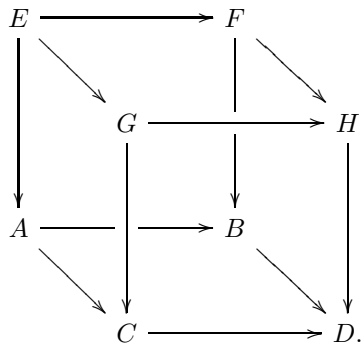
\[\diagram
      E\rrto\drto\ddto & & F\dline\drto & \\
      & G\rrto\ddto & \dto & H\ddto \\
      A\rline\drto & \rto & B\drto & \\
      & C\rrto & & D.
   \enddiagram\]
   Suppose the bottom face $A-B-C-D$ is a homotopy pushout
   and the sides $E-G-A-C$ and $E-F-A-B$ are homotopy pullbacks.
   \begin{letterlist}
      \item If the top face $E-F-G-H$ is also a homotopy pushout
            then the sides $G-H-C-D$ and $F-H-B-D$ are homotopy pullbacks.
      \item If the sides $G-H-C-D$ and $F-H-B-D$ are also homotopy
            pullbacks then the top face $E-F-G-H$ is a homotopy pushout.
   \end{letterlist}
   $\qqed$
\end{lemma}

We next set some notation. Let $X_{1}$ and $X_{2}$
be spaces, and let $1\leq j\leq 2$. Let
\(\pi_{j}:\namedright{X_{1}\times X_{2}}{}{X_{j}}\)
be the projection onto the $j^{th}$ factor and let
\(i_{j}:\namedright{X_{j}}{}{X_{1}\times X_{2}}\)
be the inclusion into the $j^{th}$ factor. Let
\(q_{j}:\namedright{X_{1}\vee X_{2}}{}{X_{j}}\)
be the pinch map onto the $j^{th}$ wedge summand. As well,
unless otherwise specified, we adopt the Milnor-Moore notation
of denoting the identity map on a space $X$ by $X$.

\begin{lemma}
   \label{porterpo}
   Let $A$, $B$ and $C$ be spaces. Define $Q$ as the homotopy pushout
   \[\diagram
        A\times B\rto^-{\ast\times B}\dto^{\pi_{1}} & C\times B\dto \\
        A\rto & Q.
     \enddiagram\]
   Then $Q\simeq (A\ast B)\vee (C\rtimes B)$.
\end{lemma}

\begin{proof}
Consider the diagram of iterated homotopy pushouts
\[\diagram
     A\times B\rto^-{\pi_{2}}\dto^{\pi_{1}} & B\rto^-{i_{2}}\dto^{\ast}
        & C\times B\dto^{s} \\
     A\rto^-{\ast} & A\ast B\rto^-{t} & \overline{Q}.
  \enddiagram\]
Here, it is well known that the left square is a
homotopy pushout, and the right homotopy pushout
defines~$\overline{Q}$. Note that
$i_{2}\circ\pi_{2}\simeq\ast\times B$. The outer rectangle in an
iterated homotopy pushout diagram is itself a homotopy pushout, so
$\overline{Q}\simeq Q$. The right pushout then shows that the
homotopy cofibre of \(\namedright{C\times B}{}{Q}\) is $\Sigma
B\vee (A\ast B)$. Thus~$t$ has a left homotopy inverse. Further,
$s\circ i_{2}\simeq\ast$ so pinching out $B$ in the right pushout
gives a homotopy cofibration \(\nameddright{C\rtimes
B}{}{Q}{r}{A\ast B}\) with $r\circ t$ homotopic to the identity
map.
\end{proof}

\begin{lemma}
   \label{gluingpo}
   Let $A$, $B$, $C$ and $D$ be spaces. Define $Q$ as the homotopy
   pushout
   \[\diagram
         A\times B\rto^-{\ast\times B}\dto^{A\times\ast}
             & C\times B\dto \\
         A\times D\rto & Q.
     \enddiagram\]
   Then $Q\simeq (A\ast B)\vee (C\rtimes B)\vee (A\ltimes D)$.
\end{lemma}

\begin{proof}
Let $Q_{1}$ be the homotopy pushout of the maps
\(\namedright{A\times D}{}{Q}\)
and
\(\namedright{A\times D}{\pi_{1}}{A}\).
Then there is a diagram of iterated homotopy pushouts
\[\diagram
     A\times B\rto^-{\ast\times B}\dto^{A\times\ast}
         & C\times B\dto \\
     A\times D\rto\dto^{\pi_{1}} & Q\dto \\
     A\rto & Q_{1}.
  \enddiagram\]
Observe that the outer rectangle is also a homotopy pushout,
so by Lemma~\ref{porterpo} we have
$Q_{1}\simeq (A\ast B)\vee (C\rtimes B)$.
Further, the outer rectangle shows that the map
\(\namedright{A}{}{Q_{1}}\)
is null homotopic. Since
\(\namedright{A\times B}{A\times\ast}{A\times D}\)
is homotopic to the composite
\(\nameddright{A\times B}{\pi_{1}}{A}{i_{1}}{A\times D}\),
there is an iterated homotopy pushout diagram
\[\diagram
     A\times B\rto^-{\ast\times B}\dto^{\pi_{1}}
        & C\times B\dto \\
     A\rto\dto^{i_{1}} & Q_{1}\dto \\
     A\times D\rto & Q.
  \enddiagram\]
Since
\(\namedright{A}{}{Q_{1}}\)
is null homotopic, we can pinch out $A$ in the lower pushout
to obtain a homotopy pushout
\[\diagram
      \ast\rto\dto & Q_{1}\dto \\
      A\ltimes D\rto & Q.
  \enddiagram\]
Hence
$Q\simeq Q_{1}\vee (A\ltimes D)\simeq
    (A\ast B)\vee (C\rtimes D)\vee (A\ltimes D)$.
\end{proof}

\begin{lemma}
   \label{shiftedpo}
   Let $A$, $B$ and $C$ be spaces. Define $Q$ as the homotopy
   pushout
   \[\diagram
        A\times (B\vee C)\rto^-{\pi_{2}}\dto^{A\times q_{2}}
           & B\vee C\dto \\
        A\times C\rto & Q.
     \enddiagram\]
   Then $Q\simeq C\vee (A\ast B)$.
\end{lemma}

\begin{proof}
First consider the homotopy pushout
\[\diagram
     B\rto\dto & B\vee C\dto^{q_{2}} \\
     \ast\rto & C.
  \enddiagram\]
In general, if $M$ is the homotopy pushout of maps
\(\namedright{X}{f}{Y}\) and \(\namedright{X}{g}{Z}\) then an easy
application of the Cube Lemma (Lemma~\ref{cube}) shows that $N\times
M$ is the homotopy pushout of \(\llnamedright{N\times X}{N\times
f}{N\times Y}\) and \(\llnamedright{N\times X}{N\times g}{N\times
Z}\). In our case, taking the product with $A$ gives a homotopy
pushout
\[\diagram
    A\times B\rto\dto^{\pi_{1}} & A\times (B\vee C)\dto^{A\times q_{2}} \\
    A\rto^-{i_{1}} & A\times C.
  \enddiagram\]
Now consider the diagram of iterated homotopy pushouts
\[\diagram
     A\times B\rto\dto^{\pi_{1}}
       & A\times (B\vee C)\rto^{\pi_{2}}\dto^{A\times q_{2}}
       & B\vee C\rto^{q_{1}}\dto & B\dto \\
     A\rto^-{i_{1}} & A\times C\rto & Q\rto & Q^{\prime}
  \enddiagram\]
where the right pushout defines $Q^{\prime}$. Because the squares
are all homotopy pushouts so is the outermost rectangle. Thus,
as the top row is homotopic to the projection $\pi_{2}$, we see that
$Q^{\prime}\simeq A\ast B$. The right pushout then implies
there is a homotopy cofibration
\(\nameddright{C}{}{Q}{}{Q^{\prime}\simeq A\ast B}\).

On the other hand, the composite
\(\nameddright{A\times B}{}{A\times (B\vee C)}{\pi_{2}}{B\vee C}\)
is homotopic to the composite
\(\nameddright{A\times B}{\pi_{2}}{B}{j_{1}}{B\vee C}\),
where $j_{1}$ is the inclusion. Thus there is an iterated
homotopy pushout diagram
\[\diagram
     A\times B\rto^-{\pi_{2}}\dto^{\pi_{1}} & B\rto^-{j_{1}}\dto
        & B\vee C\dto \\
     A\rto & A\ast B\rto & Q.
  \enddiagram\]
As $q_{1}\circ j_{1}$ is homotopic to the identity map on $B$,
the composite
\(\nameddright{A\ast B}{}{Q}{}{Q^{\prime}\simeq A\ast B}\)
is homotopic to the identity map. Hence the homotopy cofibration
\(\nameddright{C}{}{Q}{}{A\ast B}\)
splits as $Q\simeq C\vee (A\ast B)$.
\end{proof}

\begin{lemma}
   \label{poretract}
   Suppose there is a homotopy pushout
   \[\diagram
        A\times B\rto^-{f}\dto^{\ast\times B} & D\dto \\
        C\times B\rto^-{g} & E
     \enddiagram\]
   where the restriction of $f$
   to $B$ is null homotopic. Then $g$ factors through a map
   \(g^{\prime}:\namedright{C\rtimes B}{}{E}\)
   and $g^{\prime}$ has a left homotopy inverse.
\end{lemma}

\begin{proof}
As the restriction of $f$ to $B$ is null homotopic, the homotopy
commutativity of the diagram in the statement of the Lemma implies
that the restriction of $g$ to $B$ is also null homotopic. Pinching
$B$ out on the left side results in a homotopy pushout
\[\diagram
     A\rtimes B\rto^-{f^{\prime}}\dto^{\ast\rtimes B} & D\dto \\
     C\rtimes B\rto^-{g^{\prime}}\dto & E\dto \\
     Y\rdouble & Y
  \enddiagram\]
for maps $f^{\prime}$ and $g^{\prime}$. Since $\ast\rtimes B$ is
null homotopic, we have
$Y\simeq (C\rtimes B)\vee \Sigma (A\rtimes B)$, implying that
$g^{\prime}$ has a left homotopy inverse.
\end{proof}

\section{A review of homotopy actions}
\label{sec:actions}

This section is a brief reminder of some properties of homotopy
actions. Suppose there is a homotopy fibration
\[\nameddright{F}{}{E}{}{B}.\]
Let
\(\partial:\namedright{\Omega B}{}{F}\)
be the connecting map in the homotopy fibration sequence. Then
there is a canonical homotopy action
\(\theta:\namedright{F\times\Omega B}{}{F}\)
such that:
\begin{letterlist}
   \item $\theta$ restricted to $F$ is homotopic to the identity map,
   \item $\theta$ restricted to $\Omega B$ is homotopic to $\partial$, and
   \item there is a homotopy commutative diagram
         \[\diagram
              \Omega B\times\Omega B\rto^-{\mu}\dto^{\partial\times\Omega B}
                  & \Omega B\dto^{\partial} \\
              F\times\Omega B\rto^-{\theta} & F.
           \enddiagram\]
\end{letterlist}
A special case is given by the path-loop fibration
\(\nameddright{\Omega B}{}{\mathcal{P}B}{}{B}\). Here, the homotopy
action \(\theta:\namedright{\Omega B\times\Omega B}{}{\Omega B}\) is
homotopic to the loop multiplication.

Next, the homotopy action is natural for maps of homotopy fibration
sequences. If there is a homotopy fibration diagram
\[\diagram
     F\rto\dto^{f} & E\rto\dto^{g} & B\dto^{h} \\
     F^{\prime}\rto & E^{\prime}\rto & B^{\prime},
  \enddiagram\]
then there is a homotopy commutative diagram of actions
\[\diagram
     F\times\Omega B\rto^-{\theta}\dto^{f\times\Omega h} & F\dto^{f} \\
     F^{\prime}\times\Omega B^{\prime}\rto^-{\theta^{\prime}} & F^{\prime}.
  \enddiagram\]
One example of this, we will make use of, is the following.

\begin{lemma}
   \label{projaction}
   Suppose
   \(\nameddright{F}{}{E}{f}{B}\)
   is a homotopy fibration with homotopy action
   \(\theta:\namedright{F\times\Omega B}{}{F}\).
   Then the homotopy fibration
   \(\nameddright{F}{}{E\times X}{f\times X}{B\times X}\)
   has a homotopy action
   \(\theta^{\prime}:\namedright{F\times (\Omega B\times\Omega X)}{}{F}\)
   which factors as
   \[\diagram
        F\times (\Omega B\times\Omega X)\rto^-{\theta^{\prime}}
            \dto^{F\times\pi_{1}} & F\ddouble \\
        F\times\Omega B\rto^-{\theta} & F
     \enddiagram\]
   where $\pi_{1}$ is the projection.
\end{lemma}

\begin{proof}
Projecting, we obtain a homotopy pullback
\[\diagram
      F\rto\ddouble & E\times X\rto^-{f\times X}\dto^{\pi_{1}}
         & B\times X\dto^{\pi_{1}} \\
      F\rto & E\rto^-{f} & B.
  \enddiagram\]
The asserted homotopy commutative diagram now follows from the
naturality of the homotopy action.
\end{proof}

\section{A special case of the Cube Lemma}
\label{sec:cubecase}

This section describes a particular case of the Cube Lemma
which involves a homotopy action in the homotopy pushout of fibres.
Suppose there is a homotopy pushout
\[\diagram
     A\rto\dto & B\dto \\
     C\rto & D.
  \enddiagram\]
Suppose there is a space $Z$ and a map \(\namedright{D}{}{Z}\). Map
each of $A$, $B$, $C$ and $D$ into $Z$ and take homotopy fibres;
name these $E$, $F$, $G$ and $H$ respectively. Then there is a
homotopy commutative cube
   \[\diagram
      E\rrto\drto\ddto & & F\dline\drto & \\
      & G\rrto\ddto & \dto & H\ddto \\
      A\rline\drto & \rto & B\drto & \\
      & C\rrto & & D
   \enddiagram\]
in which the bottom face is a homotopy pushout and all four sides
are homotopy pullbacks. Lemma~\ref{cube} then says that the top face
is also a homotopy pushout. In practice, we will have $Z=C\times Y$
for some space $Y$, together with two additional conditions,
described in the following proposition.

\begin{proposition}
   \label{cubecase}
   Suppose there is a decomposition $Z=C\times Y$ such that:
  \begin{itemize}\item [(i)] the composite
   \(\nameddright{C}{}{D}{}{C\times Y}\)
   is homotopic to the inclusion of the first factor;
  \item [(ii)] the composite
   \(\nameddright{B}{}{D}{}{C\times Y}\)
   has a right homotopy inverse when looped.
  \end{itemize}
   Let $M$ be the homotopy fibre of the map
   \(\namedright{A}{}{C}\).
   Then:
   \begin{letterlist}
      \item $E\simeq M\times\Omega Y$ and $G\simeq\Omega Y$;
      \item the homotopy pushout of fibres becomes
            \[\diagram
                 M\times\Omega Y\rto^-{g}\dto^{\pi} & F\dto \\
                 \Omega Y\rto & H
              \enddiagram\]
            where $\pi$ is the projection and the restriction of $g$
            to $\Omega Y$ is null homotopic;
      \item the map $g$ is homotopic to the composite
            \[\llnamedddright{M\times\Omega Y}{g\vert_{M}\times\Omega Y}
                 {F\times\Omega Y}{F\times i}
                 {F\times (\Omega C\times\Omega Y)}{\theta}{F}\]
            where $g\vert_{M}$ is the restriction of $g$ to $M$,
            $i$ is the inclusion into the second factor, and $\theta$
            is the homotopy action of $\Omega C\times\Omega Y$ on $F$.
   \end{letterlist}
\end{proposition}

\begin{proof}
First consider the effect of condition~(i) on the cube, in
particular, on the face $E-G-A-C$. Let $\mathcal{P}Y$ be the path
space of $Y$. The inclusion \(\namedright{C}{}{C\times Y}\) can be
replaced up to homotopy equivalence by the product
\(\namedright{C\times \mathcal{P}Y}{}{C\times Y}\). The map
\(\namedright{A}{}{C}\) is then replaced by the product map
\(\namedright{A\times\ast}{}{C\times \mathcal{P}Y}\). Composing into
$C\times Y$ then gives a homotopy pullback
\[\diagram
      N\rto\dto & A\times\ast\rto\dto & C\times Y\ddouble \\
      \ast\times\Omega Y\rto & C\times \mathcal{P}Y\rto & C\times Y
  \enddiagram\]
which defines the space $N$. Since the maps defining the homotopy
pullback are all product maps, $N$ is homotopy equivalent to the
product $N_{1}\times N_{2}$, where $N_{1}$ is the homotopy pullback
of the maps \(\namedright{A}{}{C}\) and \(\namedright{\ast}{}{C}\),
and $N_{2}$ is the homotopy pullback of the maps
\(\namedright{\ast}{}{\mathcal{P}Y}\) and \(\namedright{\Omega
Y}{}{\mathcal{P}Y}\). That is, $N_{1}\simeq M$ and
$N_{2}\simeq\Omega Y$. Further, the map \(\namedright{N}{}{\Omega
Y}\) is homotopic to the projection \(\namedright{M\times\Omega
Y}{}{\Omega Y}\). This proves part~(a), that $E\simeq M\times\Omega
Y$ and $G\simeq\Omega Y$, and also shows in part~(b) that the map
\(\namedright{E}{}{G}\) corresponds to the projection.

Next, consider the cube face $E-F-A-B$. Observe that the connecting map
\(\namedright{\Omega C\times\Omega Y}{}{N}\)
for the fibration along the top row of the pullback defining $N$
corresponds to the product map
\(\lnamedright{\Omega C\times\Omega Y}{\delta\times\Omega Y}
   {M\times\Omega Y}\),
where $\delta$ is the connecting map in the homotopy fibration sequence
\(\namedddright{\Omega C}{\delta}{M}{}{A}{}{C}\).
Using the homotopy equivalence $E\simeq M\times\Omega Y$ we are considering
the homotopy pullback diagram
\[\diagram
     \Omega C\times\Omega Y\rto^-{\delta\times\Omega Y}\ddouble
         & M\times\Omega Y\rto\dto^{g}
         & A\rto\dto & C\times Y\ddouble \\
     \Omega C\times\Omega Y\rto^-{\gamma} & F\rto & B\rto & C\times Y
  \enddiagram\]
where $\gamma$ is the connecting map. Condition~(ii) implies
that $\gamma$ is null homotopic. The homotopy commutativity
of the left square then immediately implies that the restriction
of $g$ to $\Omega Y$ is null homotopic. This completes the proof
of part~(b).

The naturality of the homotopy action applied to the homotopy pullback
in the previous paragraph gives a homotopy commutative diagram
\[\diagram
     (M\times\Omega Y)\times (\Omega C\times\Omega Y)
          \rto^-{\theta^{\prime}}\dto^{g\times (\Omega C\times\Omega Y)}
        & M\times\Omega Y\dto^{g} \\
     F\times (\Omega C\times\Omega Y)\rto^-{\theta} & F,
  \enddiagram\]
where $\theta^{\prime}$ and $\theta$ are the respective actions.
Since the homotopy fibration sequence
\(\namedddright{\Omega C\times\Omega Y}{\delta\times\Omega Y}{M\times\Omega Y}
     {}{A\times\ast}{}{C\times Y}\)
is a product of fibration sequences, $\theta^{\prime}$ is homotopic
to the product of the actions of the individual fibrations. That is,
the homotopy fibration sequence \(\namedddright{\Omega
C}{\delta}{M}{}{A}{}{C}\) has a homotopy action
\(\theta^{\prime\prime}:\namedright{M\times\Omega C}{}{M}\), while
the homotopy fibration sequence \(\namedddright{\Omega Y}{}{\Omega
Y}{}{\mathcal{P}Y}{}{Y}\) has a homotopy action
\(\mu:\namedright{\Omega Y\times\Omega Y}{}{\Omega Y}\) given by the
loop multiplication. The map $\theta^{\prime}$ is then the composite
\[\theta^{\prime}:
    \lllnameddright{(M\times\Omega Y)\times (\Omega C\times\Omega Y)}
    {M\times T\times\Omega Y}{(M\times\Omega C)\times (\Omega Y\times\Omega Y)}
    {\theta^{\prime\prime}\times\mu}{M\times\Omega Y}\]
where $T$ is the map which interchanges factors. Precomposing
with the inclusion of factors $1$ and~$4$,
\(\namedright{M\times\Omega Y}{j\times i}
      {(M\times\Omega Y)\times (\Omega C\times\Omega Y)}\),
we have $\theta^{\prime}\circ (j\times i)$ homotopic to the identity
map. The homotopy commutative diagram of actions above then results
in a string of homotopies
\[g\simeq g\circ\theta^{\prime}\circ (j\times i)\simeq
    \theta\circ (g\times 1_{\Omega Y\times\Omega Y})\circ (j\times i)\simeq
    \theta\circ (g\vert_{M}\times i)\]
which proves part~(c).
\end{proof}

\begin{corollary}
   \label{cubecasecor}
   There is a homotopy cofibration
   \[\nameddright{M\rtimes\Omega Y}{g^{\prime}}{F}{}{H}\]
   where $g^{\prime}$ is an extension of $g$ to $M\rtimes\Omega Y$.
\end{corollary}

\begin{proof}
Consider the homotopy pushout of fibres in
Proposition~\ref{cubecase}. We know that the restriction of $g$ to
$\Omega Y$ is null homotopic. Since the projection $\pi$ has a right
inverse, the map \(\namedright{\Omega Y}{}{H}\) is also null
homotopic. Thus the factor $\Omega Y$ in the left column of the
homotopy pushout can be pinched out, resulting in a new homotopy
pushout
\[\diagram
      M\rtimes Y\rto^-{g^{\prime}}\dto & F\dto \\
      \ast\rto & H
  \enddiagram\]
which is exactly the asserted homotopy cofibration.
\end{proof}

\section{Proper coordinate subspaces of the fat wedge}
\label{sec:fatwedge}

Let $X_{1},\ldots,X_{n}$ be path-connected spaces.
In this section we investigate properties of the homotopy fibre
of the inclusion of the fat wedge $FW(1,\ldots,n)$ into the product
$X_{1}\times\cdots\times X_{n}$. Here,
\[FW(1,\ldots,n)=\{(x_{1},\ldots,x_{n})\,\vert\
    \mbox{at least one $x_{i}$ is $\ast$}\}.\]
Including the fat wedge into the product gives a homotopy fibration
\[\nameddright{F^{n}}{}{FW(1,\ldots,n)}{}
      {X_{1}\times\cdots\times X_{n}}\]
which defines the space $F^{n}$. Porter~\cite{P} showed that $F^{n}$
is homotopy equivalent to $\Omega X_{1}\ast\cdots\ast\Omega X_{n}$
by examining certain subspaces of contractible spaces. Doeraene~\cite{D}
reproduced this result in a more general setting by using the Cube Lemma.
We include a proof using the Cube Lemma for the sake of completeness.

\begin{lemma}
   \label{fatwedgefib}
   There is a homotopy equivalence
   $F^{n}\simeq\Omega X_{1}\ast\cdots\ast\Omega X_{n}$.
\end{lemma}

\begin{proof}
We induct on $n$. When $n=1$, we have $FW(1)=\ast$ and so
$F^{1}=\Omega X_{1}$. Assume $F^{n-1}\simeq\Omega
X_{1}\ast\cdots\ast\Omega X_{n-1}$. Observe that there is a
topological pushout
\[\diagram
      FW(1,\ldots,n-1)\rto^-{i}\dto & FW(1,\ldots,n-1)\times X_{n}\dto \\
      X_{1}\times\cdots\times X_{n-1}\rto & FW(1,\ldots,n)
  \enddiagram\]
where $i$ is the inclusion into the first factor. Mapping all four
corners into $X_{1}\times\cdots\times X_{n}$ and taking homotopy
fibres gives homotopy fibrations
\begin{equation}\label{1}\nameddright{F^{n}}{}{FW(1,\ldots,n)}{}{X_{1}\times\cdots\times
X_{n}}\end{equation}
\begin{equation}\label{2}\nameddright{F^{n-1}}{}{FW(1,\ldots,n-1)\times X_{n}}{}
     {X_{1}\times\cdots\times X_{n}}\end{equation}
\begin{equation}\label{3}\nameddright{\Omega X_{n}}{}{X_{1}\times\cdots\times X_{n-1}}{}
     {X_{1}\times\cdots\times X_{n}}\end{equation}
\begin{equation}\label{4}\nameddright{F^{n-1}\times\Omega X_{n}}{}{FW(1,\ldots,n-1)}{}
     {X_{1}\times\cdots\times X_{n}}.\end{equation}
Note that homotopy fibration \eqref{3} is the product of the
identity fibration \(\nameddright{\ast}{}{X_{1}\times\cdots\times
X_{n-1}}{}
     {X_{1}\times\cdots\times X_{n-1}}\)
and the path-loop fibration \(\nameddright{\Omega
X_{n}}{}{\ast}{}{X_{n}}\). This relates to both homotopy fibrations
\eqref{2} and \eqref{4}. Homotopy fibration \eqref{2} is the product
of the fibration \(\nameddright{F^{n-1}}{}{FW(1,\ldots,n-1)}{}
   {X_{1}\times\cdots\times X_{n-1}}\)
and the identity fibration above. Hence the inclusion $i$ induces a
map of fibres \(\namedright{F^{n-1}\times\Omega X_{n}}{}{F^{n-1}}\)
which is the projection onto the first factor. Homotopy fibration
\eqref{4} is the product of the fibration
\(\nameddright{F^{n-1}}{}{FW(1,\ldots,n-1)}{}
   {X_{1}\times\cdots\times X_{n-1}}\)
and the path-loop fibration. Hence the inclusion
\(\namedright{FW(1,\ldots,n-1)}{}{X_{1}\times\cdots\times X_{n-1}}\)
induces a map of fibrations
\(\namedright{F^{n-1}\times\Omega X_{n}}{}{\Omega X_{n}}\)
which is the projection onto the second factor. Collecting all this
information on the homotopy fibres, Lemma~\ref{cube} says that there
is a homotopy pushout of fibres
\[\diagram
     F^{n-1}\times\Omega X_{n}\rto^-{\pi_{1}}\dto^{\pi_{2}}
         & F^{n-1}\dto \\
     \Omega X_{n}\rto & F^{n}.
  \enddiagram\]
It is well known that in general the homotopy pushout of the projections
\(\namedright{A\times B}{}{A}\)
and
\(\namedright{A\times B}{}{B}\)
is homotopy equivalent to $A\ast B$. Thus, in our case,
$F^{n}\simeq F^{n-1}\ast\Omega X_{n}$. The inductive hypothesis on $F^{n-1}$
then implies that
$F^{n}\simeq\Omega X^{1}\ast\cdots\ast\Omega X_{n}$.
\end{proof}

For $1\leq i\leq n$, let
$X_{1}\times\cdots\times\widehat{X}_{i}\times\cdots\times X_{n}$
be the subspace of $X_{1}\times\cdots\times X_{n}$ in which the
$i^{th}$-coordinate is fixed as $\ast$. Let $FW(1,\cdots,\hat{i},\cdots n)$
be the fat wedge in
$X_{1}\times\cdots\times\widehat{X}_{i}\times\cdots\times X_{n}$.
Let $B_{i}=X_{i}\times FW(1,\ldots,\hat{i},\ldots n)$. Observe that
each $B_{i}$ is a subspace of $FW(1,\ldots,n)$ and there is a topological
pushout
\begin{equation}
  \label{Bidgrm}
  \diagram
     FW(1,\ldots,\hat{i},\ldots n)\rto\dto & B_{i}\dto \\
     X_{1}\times\cdots\times\widehat{X}_{i}\times\cdots\times X_{n}\rto
         & FW(1,\ldots,n).
  \enddiagram
\end{equation}

Consider the sequence of inclusions
\(\nameddright{B_{i}}{}{FW(1,\ldots,n)}{}{X_{1}\times\cdots\times X_{n}}\).
Using Lemma~\ref{fatwedgefib}, we obtain a homotopy pullback
\[\diagram
      F_{i}\rto^-{h_{i}}\dto & \Omega X_{1}\ast\cdots\ast\Omega X_{n}\dto \\
      B_{i}\rto\dto & FW(1,\ldots,n)\dto \\
      X_{1}\times\cdots\times X_{n}\rdouble
         & X_{1}\times\cdots\times X_{n}
  \enddiagram\]
which defines the map $h_{i}$.

\begin{lemma}
   \label{hitriv}
   The map $h_{i}$ is null homotopic.
\end{lemma}

\begin{proof}
Consider the homotopy pushout in diagram~\eqref{Bidgrm}. We wish to
apply Proposition~\ref{cubecase} with
$A=FW(1,\ldots,\hat{i},\ldots,n)$, $B=B_{i}$,
$C=X_{1}\times\cdots\times\widehat{X}_{i}\times\cdots\times X_{n}$,
$D=FW(1,\ldots,n)$, and $Z=X_{1}\times\cdots\times X_{n}$. We need
to check that the two conditions in Proposition~\ref{cubecase} hold.
Observe that $Z=C\times X_{i}$ and \(\namedright{C}{}{Z}\) is the
inclusion of the first factor so condition~(i) is satisfied. Since
$B_{i}=X_{i}\times FW(1,\ldots,\hat{i},\ldots n)$ and the map
\(\namedright{FW(1,\ldots,\hat{i},\ldots,n)}{}
   {X_{1}\times\cdots\widehat{X}_{i},\ldots\times X_{n}}\)
has a right homotopy inverse when looped, the map
\(\namedright{B_{i}}{}{X_{1}\times\cdots\times X_{n}}\) also has a
right homotopy inverse when looped, and so condition~(ii) is
satisfied. Proposition~\ref{cubecase} then says that when the four
corners of the pushout in diagram~\eqref{Bidgrm} are mapped into
$X_{1}\times\cdots\times X_{n}$ and homotopy fibres are taken, there
is a homotopy pushout of fibres
\begin{equation}
   \label{Fidgrm}
   \diagram
       (\Omega X_{1}\ast\cdots\ast\widehat{\Omega X_{i}}
            \ast\cdots\ast\Omega X_{n})
          \times\Omega X_{i}\rto^-{g}\dto^{\pi} & F_{i}\dto^{h_{i}} \\
       \Omega X_{i}\rto & \Omega X_{1}\ast\cdots\ast\Omega X_{n}
   \enddiagram
\end{equation}
where $\pi$ is the projection, the restriction of $g$ to $\Omega X_{i}$
is null homotopic, and $g$ is determined by the action of
$\Omega X_{1}\times\cdots\times\Omega X_{n}$ on $F_{i}$.

We next examine how $g$ is determined by this action.
Since $B_{i}=X_{i}\times FW(1,\ldots,\hat{i},\ldots,n)$, we
can project to obtain a homotopy pullback
\[\diagram
      F_{i}\rdouble\dto & F_{i}\dto \\
      B_{i}\rto\dto & FW(1,\ldots,\hat{i},\ldots,n)\dto \\
      X_{1}\times\cdots\times X_{n}\rto^-{\pi}
         & X_{1}\times\cdots\times\widehat{X}_{i}\times\cdots\times X_{n}.
  \enddiagram\]
Lemma~\ref{projaction} says that $g$ factors through a projection,
\[\diagram
     (\Omega X_{1}\ast\cdots\ast\widehat{\Omega X_{i}}\ast\cdots
         \ast\Omega X_{n})\times\Omega X_{i}\rto^-{g}\dto^{\pi}
         & F_{i}\ddouble \\
     (\Omega X_{1}\ast\cdots\ast\widehat{\Omega X_{i}}\ast\cdots
         \ast\Omega X_{n})\rto & F_{i}.
  \enddiagram\]
The projection of $g$ lets us define a composite
\[g^{\prime}:\nameddright{(\Omega X_{1}\ast\cdots\ast\widehat{\Omega X}_{i}
    \ast\cdots\ast\Omega X_{n})\rtimes\Omega X_{i}}{\pi}
    {\Omega X_{1}\ast\cdots\ast\widehat{\Omega X}_{i}\ast\cdots\ast
    \Omega X_{n}}{}{F_{i}}.\]
We can use $g^{\prime}$ to pinch out the factor of $\Omega X_{i}$ in
diagram~\eqref{Fidgrm} in order to obtain a homotopy cofibration
\[\nameddright{(\Omega X_{1}\ast\cdots\ast\widehat{\Omega X_{i}}
     \ast\cdots\ast\Omega X_{n})\rtimes\Omega X_{i}}{g^{\prime}}{F_{i}}{h_{i}}
     {\Omega X_{1}\ast\cdots\ast\Omega X_{n}}.\]

To simplify notation, let
$Y=\Omega X_{1}\ast\cdots\ast\widehat{\Omega X_{i}}\ast\cdots\ast\Omega X_{n}$.
Since $Y$ is a suspension,
$Y\rtimes\Omega X_{i}\simeq Y\vee (Y\wedge\Omega X_{i})$.
So $g^{\prime}$ can alternatively be described by the composite
\(\namedddright{Y\rtimes\Omega X_{i}}{\simeq}{Y\vee (Y\wedge\Omega X_{i})}
    {q}{Y}{}{F_{i}}\),
where $q$ is the pinch map. Thus $Y\wedge\Omega X_{i}$ is sent
trivially into $F_{i}$ by $g^{\prime}$ and so $\Sigma Y\wedge\Omega X_{i}$
retracts off the homotopy cofibre $\Omega X_{1}\ast\cdots\ast\Omega X_{n}$
of $g^{\prime}$. But
$\Sigma Y\wedge\Omega X_{i}\simeq\Omega X_{1}\ast\cdots\ast\Omega X_{n}$.
Thus in the homotopy cofibration sequence
\(\namedddright{Y\rtimes\Omega X_{i}}{g^{\prime}}{F_{i}}{h_{i}}
     {\Omega X_{1}\ast\cdots\ast\Omega X_{n}}{\delta}
     {\Sigma (Y\rtimes\Omega X_{i})}\),
the map $\delta$ has a left homotopy inverse and hence $h_{i}$
is null homotopic.
\end{proof}
In what follows a coordinate subspace denotes an arbitrary union of
$X_{i_1}\times\ldots\times X_{i_j}$ for some $1\leq i_1 <\ldots <
i_j\leq n$. We now use the spaces $B_{i}$ and Lemma~\ref{hitriv} to
generalise to the case of any proper coordinate subspace of
$FW(1,\ldots,n)$.

\begin{proposition}
   \label{htriv}
   Suppose $A$ is a proper coordinate subspace of $FW(1,\ldots,n)$.
   Include $A$ into $FW$ and then include into
   $X_{1}\times\cdots\times X_{n}$ to obtain a homotopy pullback
   \[\diagram
        F\rto^-{h}\dto & \Omega X_{1}\ast\cdots\ast\Omega X_{n}\dto \\
        A\rto\dto & FW(1,\ldots,n)\dto \\
        X_{1}\times\cdots\times X_{n}\rdouble
           & X_{1}\times\cdots\times X_{n}
     \enddiagram\]
   which defines the map $h$. Then $h$ is null homotopic.
\end{proposition}

\begin{proof}
First observe that the inclusion of $A$ into $FW(1,\ldots,n)$ factors
through $B_{i}$ for some $i$. This statement really just
follows from the definitions. In terms of coordinates,
\[B_{i}=\{(x_{1},\ldots,x_{m})\, \vert\ \mbox{at least one of
     $x_{1},\ldots,\hat{x}_{i},\ldots,x_{m}$ is $\ast$}\}.\]
If the inclusion of $A$ into $FW(1,\ldots,n)$ does not factor
through $B_{i}$, then $A$ must contain a sequence of the form
$(x_{1},\ldots,x_{n})$ in which each of
$x_{1},\ldots,\hat{x}_{i},\ldots,x_{n}$ is not~$\ast$. Since $A$ is
a coordinate subspace, every sequence of this form must be in $A$.
(Note that $A$ is a subspace of $FW(1,\ldots,n)$ so this forces
$x_{i}$ to be $\ast$ in each such sequence.) If this is true for
$1\leq i\leq n$, then all of $FW(1,\ldots,n)$ is contained in $A$,
contradicting the hypothesis that $A$ is a proper coordinate
subspace of~$FW(1,\ldots,n)$.

The factorization of
\(\namedright{A}{}{FW(1,\ldots,n)}\)
through $B_{i}$ results in a diagram of iterated homotopy pullbacks
\[\diagram
      F\rto\dto & F_{i}\rto^-{h_{i}}\dto
         & \Omega X_{1}\ast\cdots\ast\Omega X_{n}\dto \\
      A\rto\dto & B_{i}\rto\dto & FW(1,\ldots,n)\dto \\
      X_{1}\times\cdots\times X_{n}\rdouble
         & X_{1}\times\cdots\times X_{n}\rdouble
         & X_{1}\times\cdots\times X_{n}.
  \enddiagram\]
The outer rectangle is the homotopy pullback defining $h$,
so $h$ factors through $h_{i}$. But $h_{i}$ is null
homotopic by Lemma~\ref{hitriv}, and so $h$ is null homotopic.
\end{proof}

\section{Homotopy fibres associated to regular sequences}
\label{sec:fibreg}

Let $X_{1},\ldots,X_{n}$ be path-connected spaces. Let $A$ and $B$
be two coordinate subspaces of $X_{1}\times\cdots\times X_{n}$,
where $B\subseteq A$. Let $F_{A}$ and $F_{B}$ be the homotopy fibres
of the inclusions of $A$ and $B$ respectively into
$X_{1}\times\cdots\times X_{n}$. Observe that there is a map of
fibres \(\namedright{F_{B}}{}{F_{A}}\). The purpose of this section
is to consider the homotopy types of $F_{A}$ and $F_{B}$ and how
these are related by the map of fibres. In general, not much could
be expected to be said. We show that if $A$ is built up from $B$ by
what we call a regular sequence, and if the homotopy type of $F_{B}$
is of a certain description, then the homotopy type of $F_{A}$ is of
the same description and there is control over the map of fibres.
All this is made concrete in Theorem~\ref{construction} and
Proposition~\ref{F0toFl}.

We begin by defining what is meant by a regular sequence.
Let $\{i_{1},\ldots,i_{m}\}$ be a subset of $\{1,\ldots,n\}$,
where $i_{1}<\cdots <i_{m}$. Let $\{j_{1},\ldots,j_{n-m}\}$
be the complement of $\{i_{1},\ldots,i_{m}\}$ in $\{1,\ldots,n\}$,
where $j_{1}<\cdots <j_{n-m}$. Let $FW(i_{1},\dots,i_{m})$ be
the fat wedge in $X_{i_{1}}\times\cdots\times X_{i_{m}}$. Let
$A_{0}$ and $A$ be coordinate subspaces of $X_{1}\times\cdots\times X_{n}$
such that $X_{1}\vee\cdots\vee X_{n}\subseteq A_{0}$ and $A_{0}\subseteq A$.
Then $A$ can be built up iteratively from $A_{0}$ by a sequence of
topological pushouts
\begin{equation}
  \label{Akpo}
  \diagram
      FW(i_{1},\ldots,i_{m})\rto\dto & A_{k-1}\dto \\
      X_{i_{1}}\times\cdots\times X_{i_{m}}\rto & A_{k}
  \enddiagram
\end{equation}
where $1\leq k\leq l$, and $A_{l}=A$. There may be many choices of
sequences of pushouts which realise $A$ in this way. A particular
type of sequence, if it exists, is well suited to identifying the
homotopy fibre of the inclusion
\(\namedright{A}{}{X_{1}\times\cdots\times X_{n}}\).

\begin{definition}
Let $X_{1},\ldots,X_{n}$ be path-connected spaces.
A coordinate subspace $A$ of $X_{1}\times\cdots\times X_{n}$
is \emph{regular} if the sequence
\[A_{0}\subseteq A_{1}\subseteq\cdots\subseteq A_{l}=A\]
has the following property for each $1\leq k\leq l$. Let
$\{s_{1},\ldots,s_{r}\}$ be the largest subset of
$\{j_{1},\ldots,j_{n-m}\}$ for which $A_{k-1}$ can be written as a
product $A_{k-1}=A^{\prime}_{k-1}\times X_{s_{1}}\times\cdots\times
X_{s_{r}}$ (permuting the coordinates if necessary). Then there is a
topological pushout
\[\diagram
     M_{k-1}\rto\dto & N_{k-1}\dto \\
     FW(i_{1},\ldots,i_{m})\rto & A_{k-1}
  \enddiagram\]
where $M_{k-1}$ is a proper coordinate subspace of $FW(i_{1},\ldots,i_{m})$.
\end{definition}

The definition of a regular sequence may seem on first reading to
be a bit mystifying, but it arises naturally when considering
coordinate subspaces associated to shifted complexes. It might
be useful at this point to briefly skip ahead to Examples~\ref{regex1}
and~\ref{regex2} in order to see the connection.

To go along with the definition, we establish some notation.
Let $\{t_{1},\ldots,t_{n-m-r}\}$ be the complement of
$\{s_{1},\ldots,s_{r}\}$ in $\{j_{1},\ldots,j_{n-m}\}$.
Let $S=X_{s_{1}}\times\cdots\times X_{s_{r}}$ and
$T=X_{t_{1}}\times\cdots\times X_{t_{n-m-r}}$, so
$S\times T=X_{j_{1}}\times\cdots\times X_{j_{n-m}}$ and
$A_{k-1}=A^{\prime}_{k-1}\times S$.

For $0\leq k\leq l$, let $F_{k}$ be the homotopy fibre of the
inclusion \(\namedright{A_{k}}{}{X_{1}\times\cdots\times X_{n}}\).
Observe that if $A_{k-1}=A_{k-1}^{\prime}\times S$ then there is a
diagram of iterated homotopy pullbacks
\[\diagram
     F_{k-1}\rdouble\dto & F_{k-1}\rdouble\dto & F_{k-1}\dto \\
     A_{k-1}^{\prime}\rto^-{i}\dto & A_{k-1}\rto^-{\pi}\dto
        & A_{k-1}^{\prime}\dto \\
     X_{i_{1}}\times\cdots\times X_{i_{m}}\times T\rto^-{i}
     & X_{i_{1}}\times\cdots\times X_{i_{m}}\times S\times T\rto^-{\pi}
     & X_{i_{1}}\times\cdots\times X_{i_{m}}\times T
  \enddiagram\]
where $i$ and $\pi$ are the inclusion and projection respectively.

In Theorem~\ref{construction} we make the seemingly odd assumption
that the fibre $F_{0}$ is a co-$H$ space. However, in the context
of coordinate subspace arrangements, this condition arises naturally,
as we are trying to show that certain homotopy fibres (labelled
$\mathcal{Z}_{K}$) are homotopy equivalent to wedges of spheres for
appropriate simplicial complexes $K$, in which case the fibres
$\mathcal{Z}_{K}$ are co-$H$ spaces.

\begin{theorem}
   \label{construction}
   Suppose there is a regular sequence of coordinate subspaces
   \[A_{0}\subseteq A_{1}\subseteq\cdots\subseteq A_{l}=A.\]
   Assume that $F_{0}$ is a co-$H$ space. Then the following hold:
   \begin{letterlist}
      \item for $1\leq k\leq l$, there is a homotopy cofibration
            \[\nameddright{(\Omega X_{i_{m}}\ast\cdots\ast\Omega X_{i_{m}})
              \rtimes (\Omega X_{j_{1}}\times\cdots\times\Omega X_{j_{n-m}})}
              {}{F_{k-1}}{}{F_{k}}\]
            and a homotopy decomposition
            \[(\Omega X_{i_{m}}\ast\cdots\ast\Omega X_{i_{m}})
              \rtimes (\Omega X_{j_{1}}\times\cdots\times\Omega X_{j_{n-m}})
              \simeq C_{k-1}\vee D_{k-1}\]
            where $C_{k-1}$ maps trivially into $F_{k-1}$ and $D_{k-1}$
            retracts off $F_{k-1}$;
      \item there is a homotopy decomposition
            $F_{k-1}\simeq D_{k-1}\vee E_{k-1}$ for some space
            $E_{k-1}$;
      \item $F_{k}$ is a co-$H$ space and there is a homotopy decomposition
            $F_{k}\simeq\Sigma C_{k-1}\vee E_{k-1}$.
   \end{letterlist}
\end{theorem}

\begin{proof}
As the proof of part~(a) is lengthy, we begin by assuming that
part~(a) has been proved and show that parts~(b) and~(c) hold. With
$F_{0}$ as the base case, we inductively assume that $F_{k-1}$ is a
co-$H$ space. Let $E_{k-1}$ be the homotopy cofibre of
\(\namedright{D_{k-1}}{}{F_{k-1}}\). By part~(a), this map has a
left homotopy inverse \(\namedright{F_{k-1}}{}{D_{k-1}}\). Since
$F_{k-1}$ is a co-$H$ space, we can add to obtain a composite
\(\nameddright{F_{k-1}}{}{F_{k-1}\vee F_{k-1}}{}{D_{k-1}\vee
E_{k-1}}\) which is a homotopy equivalence. This proves part~(b).
Next, including $D_{k-1}$ into $C_{k-1}\vee D_{k-1}$ we obtain a
homotopy pushout
\[\diagram
      D_{k-1}\rto\ddouble & C_{k-1}\vee D_{k-1}\rto\dto & C_{k-1}\dto \\
      D_{k-1}\rto & F_{k-1}\rto\dto & E_{k-1}\dto \\
      & F_{k}\rdouble & F_{k}.
  \enddiagram\]
By part~(a) the map \(\namedright{C_{k-1}}{}{F_{k-1}}\) is null
homotopic, so in the pushout the map
\(\namedright{C_{k-1}}{}{E_{k-1}}\) is also null homotopic. Hence
$F_{k}\simeq\Sigma C_{k-1}\vee E_{k-1}$. Finally, $E_{k-1}$ is a
retract of $F_{k-1}$ which has been inductively assumed to be a
co-$H$ space, so $E_{k-1}$ is also a co-$H$ space. Thus $F_{k}$ is a
wedge of two co-$H$ spaces and so is itself a co-$H$ space.

We now prove part~(a).

\noindent \textbf{Step 1. Setting up:} Consider the pushout in
diagram~\eqref{Akpo}. We apply Proposition~\ref{cubecase} with
$A=FW(i_{1},\ldots,i_{m})$, $B=A_{k-1}$,
$C=X_{i_{1}}\times\cdots\times X_{i_{m}}$, $D=A_{k}$ and
$Z=X_{1}\times\cdots\times X_{n}$. We need to check that
conditions~(i) and~(ii) of Proposition~\ref{cubecase} hold. Observe
that $Z=C\times (X_{j_{1}}\times\cdots\times X_{i_{n-m}})$ and
\(\namedright{C}{}{Z}\) is the inclusion of the first factor, so
condition~(i) is satisfied. Since $X_{1}\vee\cdots\vee X_{n}$ is a
subspace of $A_{k-1}$ and the inclusion
\(\namedright{X_{1}\vee\cdots\vee X_{n}}{}{X_{1}\times\cdots\times
X_{n}}\) has a right homotopy inverse when looped, the inclusion
\(\namedright{A_{k-1}}{}{X_{1}\times\cdots\times X_{n}}\) also has a
right homotopy inverse when looped, and so condition~(ii) is also
satisfied. Proposition~\ref{cubecase} insures that when the four
corners of the pushout in diagram~\eqref{Akpo} are mapped into
$X_{1}\times\cdots\times X_{n}$ and homotopy fibres are taken, there
is a homotopy pushout of fibres
\begin{equation}
  \label{Fkpo}
  \diagram
      (\Omega X_{i_{1}}\ast\cdots\ast\Omega X_{i_{m}})\times
           (\Omega X_{j_{1}}\times\cdots\times\Omega X_{j_{n-m}})
           \rto^-{g}\dto^{\pi}
         & F_{k-1}\dto \\
      \Omega X_{j_{1}}\times\cdots\Omega X_{j_{n-m}}\rto & F_{k}
  \enddiagram
\end{equation}
where $\pi$ is the projection, the restriction of $g$ to $\Omega
X_{j_{1}}\times\cdots\times\Omega X_{j_{n-m}}$ is null homotopic,
and $g$ is determined by the action of $\Omega
X_{1}\times\cdots\times\Omega X_{n}$ on $F_{k-1}$. As the
restriction of $g$ to $\Omega X_{j_{1}}\times\cdots\times\Omega
X_{j_{n-m}}$ is null homotopic, we can pinch out this factor in
diagram~\eqref{Fkpo} and, as in Corollary~\ref{cubecasecor}, obtain
a homotopy cofibration
\begin{equation}
  \label{Fkcofib}
  \nameddright{(\Omega X_{i_{1}}\ast\cdots\ast\Omega X_{i_{t}})\rtimes
     (\Omega X_{j_{1}}\times\cdots\times\Omega X_{j_{n-t}})}{g^{\prime}}
     {F_{k-1}}{}{F_{k}}.
\end{equation}
where $g^{\prime}$ is an extension of $g$ to the half-smash.

\noindent
\textbf{Step 2. The summand $C_{k-1}$:}
The decomposition $A_{k-1}=A_{k-1}^{\prime}\times S$ implies that
there is a homotopy pullback
\[\diagram
     F_{k-1}\rdouble\dto & F_{k-1}\dto \\
     A_{k-1}\rto^-{\pi}\dto & A_{k-1}^{\prime}\dto \\
     X_{1}\times\cdots\times X_{n}\rto^-{\pi}
       & (X_{i_{1}}\times\cdots\times X_{i_{m}})\times T
  \enddiagram\]
where $\pi$ is the projection. Lemma~\ref{projaction} says that the
map $g$ in diagram~\eqref{Fkpo} factors through a projection,
\[\diagram
      (\Omega X_{i_{1}}\ast\cdots\ast\Omega X_{i_{m}})\times
          (\Omega S\times\Omega T)\rto^-{g}\dto^{1\times\pi}
           & F_{k-1}\ddouble \\
      (\Omega X_{i_{1}}\ast\cdots\ast\Omega X_{i_{m}})\times\Omega T
           \rto^-{\tilde{g}} & F_{k}
  \enddiagram\]
where $\tilde{g}$ is the restriction of $g$ to
$(\Omega X_{i_{1}}\ast\cdots\ast\Omega X_{i_{m}})\times\Omega T$.

Let $Y=\Omega X_{i_{1}}\ast\cdots\ast\Omega X_{i_{m}}$. Since the
restriction of $g$ to $(\Omega S\times\Omega T)$ is null homotopic,
the restriction of $\tilde{g}$ to $\Omega T$ is null homotopic. Thus
$\tilde{g}$ factors through $Y\rtimes\Omega T$. It was only
necessary to choose some extension $g^{\prime}$ of $g$ to the
half-smash in~\eqref{Fkcofib} in order to obtain the homotopy
cofibration, so we could have taken $g^{\prime}$ to be the composite
\(\nameddright{Y\rtimes (\Omega S\times\Omega T)}{1\times\pi}
   {Y\rtimes\Omega T}{}{F_{k-1}}\).
Since $Y$ is a suspension,
\[Y\rtimes (\Omega S\times\Omega T)\simeq
   Y\vee (Y\wedge\Omega S)\vee (Y\wedge\Omega T)\vee
   (Y\wedge\Omega S\wedge\Omega T).\]
Let
$C_{k-1}= (Y\wedge\Omega S)\vee (Y\wedge\Omega S\wedge\Omega T)$
and let
$D_{k-1}=Y\vee (Y\wedge\Omega T)$.
Then $g^{\prime}$ can alternatively be described by the composite
\(\nameddright{Y\rtimes (\Omega S\times\Omega T)}{\simeq}{C_{k-1}\vee D_{k-1}}
    {q}{D_{k-1}}\),
where $q$ is the pinch map. Thus $C_{k-1}$ is sent trivially
into $F_{k-1}$ by $g^{\prime}$, as asserted.

\noindent \textbf{Step 3. The summand $D_{k-1}$:} It remains to show
that $D_{k-1}=Y\vee (Y\wedge\Omega T)$ is a retract of $F_{k-1}$.
Again, we consider $A_{k-1}=A_{k-1}^{\prime}\times S$ where
$S=X_{s_{1}}\times\cdots\times X_{s_{r}}$. Observe that
$\{i_{1},\dots,i_{m}\}$ and $\{s_{1},\ldots,s_{r}\}$ are disjoint
sets in $\{1,\ldots,n\}$ so the inclusion
\(\namedright{FW(i_{1},\ldots,i_{m})}{}{A_{k-1}}\) of
diagram~\eqref{Akpo} factors as a composite
\(\nameddright{FW(i_{1},\ldots,i_{m})}{}{A_{k-1}^{\prime}}{}{A_{k-1}}\).
Define the space $A_{k}^{\prime\prime}$ as the topological pushout
\[\diagram
     FW(i_{1},\ldots,i_{m})\rto\dto & A_{k-1}^{\prime}\dto \\
     X_{i_{1}}\times\cdots\times X_{i_{m}}\rto & A_{k}^{\prime\prime}.
  \enddiagram\]
Since $A$ is regular, there is a topological pushout
\begin{equation}
  \label{Akprimepo}
  \diagram
     M_{k-1}\rto\dto & N_{k-1}\dto \\
     FW(i_{1},\ldots,i_{m})\rto & A_{k-1}^{\prime}
  \enddiagram
\end{equation}
where $M_{k-1}$ is a proper coordinate subspace of
$FW(i_{1},\ldots,i_{m})$. Note that all the spaces in
diagram~\eqref{Akprimepo} are coordinate subspaces of
$X_{i_{1}}\times\cdots\times X_{i_{m}}\times T$. We intend to map
the four corners of the pushout into $X_{i_{1}}\times\cdots\times
X_{i_{m}}\times T$, take homotopy fibres, and apply the Cube Lemma.
Before doing so we identify the homotopy fibres. Let $F_{M}$ be the
homotopy fibre of the inclusion
\(\namedright{M_{k-1}}{}{X_{i_{1}}\times\cdots\times X_{i_{m}}}\).
By Lemma~\ref{fatwedgefib}, the homotopy fibre of the inclusion
\(\namedright{FW(i_{1},\ldots,i_{m})}{}
    {X_{i_{1}}\times\cdots\times X_{i_{m}}}\)
is homotopy equivalent to $\Omega X_{i_{1}}\ast\cdots\ast\Omega X_{i_{m}}$.
Including $X_{i_{1}}\times\cdots\times X_{i_{m}}$ into
$X_{i_{1}}\times\cdots\times X_{i_{m}}\times T$
we obtain a homotopy pullback
\[\diagram
     F_{M}\times\Omega T\rto\dto^{h\times\Omega T} & M_{k-1}\rto\dto
        & X_{i_{1}}\times\cdots\times X_{i_{m}}\times T\ddouble \\
     (\Omega X_{i_{1}}\ast\cdots\ast\Omega X_{i_{m}})\times\Omega T\rto
        & FW(i_{1},\ldots,i_{m})\rto
        & X_{i_{1}}\times\cdots\times X_{i_{m}}\times T
  \enddiagram\]
for some map $h$. Let $F_{N}$ be the homotopy fibre of the inclusion
\(\namedright{N_{k-1}}{}{X_{i_{1}}\times\cdots\times X_{i_{m}}\times T}\).
The definition of a regular sequence includes the hypothesis that
$X_{1}\vee\cdots\vee X_{n}\subseteq A_{0}$, and so
$X_{1}\vee\cdots\vee X_{n}\subseteq A_{k-1}$. Having projected away
from coordinates $s_{1},\ldots, s_{l}$, we have
$X_{i_{1}}\vee\cdots\vee X_{i_{m}}\vee X_{j_{1}}\vee\cdots\vee X_{j_{t}}
    \subseteq A^{\prime}_{k-1}$.
As diagram~\eqref{Akprimepo} is a homotopy pushout and
$FW(i_{1},\ldots,i_{m})$ intersects $X_{j_{1}}\vee\cdots\vee
X_{j_{t}}$ at a point, we must have $X_{j_{1}}\vee\cdots\vee
X_{j_{t}}\subseteq N_{k-1}$. Thus $\Omega T=\Omega
X_{j_{1}}\times\cdots\times\Omega X_{j_{t}}$ retracts off $\Omega
N_{k-1}$. Therefore, in the homotopy pullback
\[\diagram
      F_{M}\times\Omega T\rto\dto^{\overline{f}} & M_{k-1}\rto\dto
         & X_{i_{1}}\times\cdots\times X_{i_{m}}\times T\ddouble \\
      F_{N}\rto & N_{k-1}\rto
         & X_{i_{1}}\times\cdots\times X_{i_{m}}\times T
  \enddiagram\]
(the pullback defines the map $\overline{f}$) the restriction
of $\overline{f}$ to $\Omega T$ is null homotopic. Now recall
from Step~$2$ that the homotopy fibre of the inclusion
\(\namedright{A_{k-1}^{\prime}}{}
    {X_{i_{1}}\times\cdots\times X_{i_{m}}\times T}\)
is homotopy equivalent to~$F_{k-1}$. Thus, when the four corners of
the pushout in diagram~\eqref{Akprimepo} are mapped into
$X_{i_{1}}\times\cdots\times X_{i_{m}}\times T$ and homotopy fibres
are taken, Lemma~\ref{cube} implies that there is a homotopy pushout
of fibres
\begin{equation}
  \label{hpo}
  \diagram
     F_{M}\times\Omega T\rto^-{\overline{f}}\dto^{h\times\Omega T}
        & F_{N}\dto \\
     (\Omega X_{i_{1}}\ast\cdots\ast\Omega X_{i_{m}})\times\Omega T
        \rto^-{\overline{g}} & F_{k-1}
  \enddiagram
\end{equation}
for some map $\overline{g}$. We can identify $\overline{g}$: it is
the restriction of the map $g$ in diagram~\eqref{Fkpo} to $(\Omega
X_{i_{1}}\ast\cdots\ast\Omega X_{i_{m}})\times\Omega T$. This is
because, as in the proof of Proposition~\ref{cubecase}~(c), the map
$\overline{g}$ is determined by the action of $\Omega
X_{i_{1}}\times\cdots\times\Omega X_{i_{m}}\times\Omega T$ on
$F_{k-1}$. But the pullback
\[\diagram
      F_{k-1}\rto\ddouble & A_{k-1}^{\prime}\rto\dto
         & X_{i_{1}}\times\cdots\times X_{i_{m}}\times T\dto \\
      F_{k-1}\rto & A_{k}\rto
         & X_{i_{1}}\times\cdots\times X_{i_{m}}\times S\times T
  \enddiagram\]
obtained from including $A_{k-1}^{\prime}$ into
$A_{k-1}=A_{k-1}^{\prime}\times S$ implies that the action of
$\Omega X_{i_{1}}\times\cdots\times\Omega X_{i_{m}}\times\Omega T$
on $F_{k-1}$ is the restriction of the action of
$\Omega X_{i_{1}}\times\cdots\times\Omega X_{i_{m}}\times
   \Omega T\times\Omega S$
on $F_{k-1}$, that is, the action of
$\Omega X_{1}\times\cdots\times\Omega X_{n}$
on $F_{k-1}$, and the latter action determines $g$.
Consequently, the factorization of $g$ through $g^{\prime}$
implies that the restriction $\overline{g}$ factors as a composite
(using the notation from Step~$2$)
\(\gamma:\nameddright{Y\rtimes\Omega T}{}{Y\rtimes(\Omega S\times\Omega T)}
    {g^{\prime}}{F_{k-1}}\).
Note that $D_{k-1}$ was defined as
$Y\vee (Y\wedge\Omega T)\simeq Y\rtimes\Omega T$,
and we are trying to prove precisely that $\gamma$ has a left
homotopy inverse.

Consider diagram~\eqref{hpo}. Since $M_{k-1}$ is a proper coordinate
subspace of $FW(i_{1},\ldots,i_{m})$, Proposition~\ref{htriv}
implies that $h$ is null homotopic. Thus $h\times\Omega T$ is
homotopic to $\ast\times\Omega T$. We have seen that the restriction
of $\overline{f}$ to $\Omega T$ is null homotopic.
Lemma~\ref{poretract} now applies, and shows that $\gamma$ has a
left homotopy inverse.
\end{proof}

We now condense some of the information coming out of
Theorem~\ref{construction} by concentrating on how the fibre
$F_{0}$ of the starting point $A_{0}$ of the regular sequence
relates to the fibre $F_{l}$ of the ending point $A_{l}$ of
the sequence. Let $\theta$ be the composite
\[\theta:\namedright{F_{0}}{}{F_{1}}\longrightarrow\cdots
    \longrightarrow F_{l}.\]
In particular, we want to know how the homotopy type of $F_{0}$ influences
that of $F_{l}$. This requires a suitable hypothesis on the homotopy
type of $F_{0}$ to get going. We now define a class of spaces which will
do the job.

\begin{definition}
\label{Gdefn}
Let $\mathcal{G}_{1}^{n}$ be the collection of spaces $F$ which
are homotopy equivalent to a wedge of summands of the form
$\Sigma\Omega X_{i_{1}}\ast\cdots\ast\Omega X_{i_{m}}$, where
$1\leq i_{1}<\cdots <i_{m}\leq n$.
\end{definition}

Consider Definition~\ref{Gdefn} in the case of primary interest, when
$X_{i}=\mathbb{C}P^{\infty}$ for $1\leq i\leq n$. Then
$\Omega X_{i}\simeq S^{1}$ and so
$\Sigma\Omega X_{i_{1}}\ast\cdots\ast\Omega X_{i_{m}}\simeq S^{m+1}$,
in which case $F$ is homotopy equivalent to a wedge of spheres.
As spaces which are homotopy equivalent to wedges of spheres will
appear repeatedly, it will be convenient to introduce an abbreviated
way of saying this.

\begin{definition}
Let $\mathcal{W}$ be the collection of spaces $F$ which are homotopy
equivalent to a wedge of spheres.
\end{definition}

\begin{proposition}
   \label{F0toFl}
   Assume the hypotheses of Theorem~\ref{construction}. Suppose in
   addition that (for all path-connected spaces $X_{1},\ldots,X_{n}$)
   we have $F_{0}\in\mathcal{G}_{1}^{n}$. Consider the map of fibres
   \(\theta:\namedright{F_{0}}{}{F_{l}}\). The following hold:
   \begin{letterlist}
      \item $F_{l}\in\mathcal{G}_{1}^{n}$, and
      \item there is a homotopy decomposition
            $F_{0}\simeq F_{0}^{1}\vee F_{0}^{2}$
            where $F_{0}^{1},F_{0}^{2}\in\mathcal{G}_{1}^{n}$, the restriction
            of $\theta$ to $F_{0}^{1}$ is null homotopic, and the restriction
            of $\theta$ to $F_{0}^{2}$ has a right homotopy inverse.
   \end{letterlist}
\end{proposition}

\begin{proof}
Theorem~\ref{construction} gives that for $1\leq k\leq l$, there is
a homotopy cofibration
\[\nameddright{(\Omega X_{i_{m}}\ast\cdots\ast\Omega X_{i_{m}})
              \rtimes (\Omega X_{j_{1}}\times\cdots\times\Omega X_{j_{n-m}})}
              {f_{k}}{F_{k-1}}{g_{k}}{F_{k}}\]
and there are homotopy decompositions
\[(\Omega X_{i_{m}}\ast\cdots\ast\Omega X_{i_{m}})
     \rtimes (\Omega X_{j_{1}}\times\cdots\times\Omega X_{j_{n-m}})
     \simeq C_{k-1}\vee D_{k-1},\]
\[F_{k-1}\simeq D_{k-1}\vee E_{k-1},\]
\[F_{k}\simeq\Sigma C_{k-1}\vee E_{k-1}\]
where the restriction of $f_{k}$ to $C_{k-1}$ is null homotopic and
the restriction of $f_{k}$ to $D_{k-1}$ has a right homotopy inverse.
This implies that the restriction of $g_{k}$ to $D_{k-1}$ is null homotopic
and the restriction of $g_{k}$ to $E_{k-1}$ has a right homotopy inverse.

First observe that, in general, there are homotopy decompositions
$(\Sigma A)\rtimes B\simeq\Sigma A\vee (\Sigma A\wedge B)$ and
$\Sigma (A\times B)\simeq\Sigma A\vee\Sigma B\vee (\Sigma A\wedge B)$.
Using the first decomposition and iterating on the second, we see that
\[(\Omega X_{i_{m}}\ast\cdots\ast\Omega X_{i_{m}})
     \rtimes (\Omega X_{j_{1}}\times\cdots\times\Omega X_{j_{n-m}})
     \in\mathcal{G}_{1}^{n}.\]
The fact that Theorem~\ref{construction} holds for all path-connected
spaces $X_{1},\ldots,X_{n}$ means that the decompositions are
independent of the particular choices of those spaces. This
lets us make an advantageous choice of $X_{1},\ldots,X_{n}$, observe
how the decompositions behave in this special case, and then infer the
general decompositions.

The advantageous choice is to take $X_{i}=\mathbb{C}P^{\infty}$
for $1\leq i\leq n$. Then $\Omega X_{i}\simeq S^{1}$ and
\[(\Omega X_{i_{m}}\ast\cdots\ast\Omega X_{i_{m}})
      \rtimes (\Omega X_{j_{1}}\times\cdots\times\Omega X_{j_{n-m}})
      \in\mathcal{W}.\]
Thus $C_{k-1},D_{k-1}\in\mathcal{W}$. The hypothesis on $F_{0}$
implies that in this case $F_{0}\in\mathcal{W}$. Thus the homotopy
equivalence $F_{0}\simeq D_{0}\vee E_{0}$ implies that
$E_{0}\in\mathcal{W}$. Hence $F_{1}\in\mathcal{W}$. Inductively, we
see that $C_{k-1},D_{k-1},E_{k-1},F_{k}\in\mathcal{W}$ for all
$1\leq k\leq l$. In particular, $F_{l}\in\mathcal{W}$. We next
describe the decomposition $F_{0}\simeq F_{0}^{1}\vee F_{0}^{2}$.
The decomposition $F_{k}\simeq\Sigma C_{k-1}\vee E_{k-1}$ gives a
retraction of $E_{k-1}$ off $F_{k}$. Consider how this relates to
the decomposition $F_{k}\simeq D_{k}\vee E_{k}$. Since
$D_{k},E_{k}\in\mathcal{W}$, we can choose subwedges of spheres
$E_{D_{k}},E_{E_{k}}$ of $D_{k},E_{k}$ respectively such that
$E_{k-1}\simeq E_{D_{k}}\vee E_{E_{k}}$. Let $F_{0}^{1}=D_{1}\vee
E_{D_{1}}\vee\cdots\vee E_{D_{l-1}}$, and let
$F_{0}^{2}=E_{E_{l-1}}$. Then $F_{0}^{1}\vee F_{0}^{2}\simeq F_{0}$.
The condition that $g_{k}$ is null homotopic when restricted to
$D_{k}$ then implies that it is null homotopic when restricted to
$E_{D_{k}}$, and so collectively we see that the restriction of
$\theta$ to $F_{0}^{1}$ is null homotopic. The condition that
$g_{l-1}$ has a right homotopy inverse when restricted to $E_{l-1}$
implies that the restriction of $\theta$ to $F_{0}^{2}=E_{E_{l-1}}$
has a right homotopy inverse.

Now consider the general case. Observe that by keeping track of the
indices $i_{s}$ and $j_{t}$ on each copy of $\Omega X_{i_{s}}\simeq S^{1}$
and $\Omega X_{j_{t}}\simeq S^{1}$ in the special case, we can discern
which wedge summands of
$(\Omega X_{i_{m}}\ast\cdots\ast\Omega X_{i_{m}})
      \rtimes (\Omega X_{j_{1}}\times\cdots\times\Omega X_{j_{n-m}})$
are in $C_{k-1}$ and which are in $D_{k-1}$. In particular, we see
that $C_{k-1},D_{k-1}\in\mathcal{G}_{1}^{n}$ for each $1\leq k\leq l$.
The same index bookkeeping on the successive decompositions in the
special case then implies that $E_{k-1},F_{k}\in\mathcal{G}_{1}^{n}$ for
$1\leq k\leq l$ -- in particular, $F_{l}\in\mathcal{G}_{1}^{n}$, proving
part~(a) -- and there is a decomposition
$F_{0}\simeq F_{0}^{1}\vee F_{0}^{2}$ such that
$F_{0}^{1},F_{0}^{2}\in\mathcal{G}_{1}^{n}$, the restriction
of $\theta$ to $F_{0}^{1}$ is null homotopic and the restriction
of $\theta$ to $F_{0}^{2}$ has a right homotopy inverse, which
proves part~(b).
\end{proof}

\section{The existence of regular sequences}
\label{sec:regseqs}

In this section we give a general set of conditions which guarantees
the existence of regular sequences. In Examples~\ref{regex1}
and~\ref{regex2} we then give particular instances which will be
used later in Section~\ref{sec:ZKtype}. The set of conditions is
phrased in terms of shifted complexes. Recall that a simplicial
complex $K$ is \emph{shifted} if there is an ordering on its set of
vertices such that whenever $\sigma\in K$ and $v^{\prime}<v$, then
$(\sigma-v)\cup v^{\prime}\in K$. Two additional definitions we
need are the following.
Let $K$ be a simplicial complex. The \textit{link} and the \textit{star}
of a simplex $\sigma \in K$ are the subcomplexes
\[\mathrm{link}_{K}\sigma =\{\tau \in K\,\vert\ \sigma \cup \tau \in K,\sigma \cap
\tau =\emptyset \};\]
\[\mathrm{star}_{K}\sigma =\{\tau \in K\,\vert\ \sigma \cup \tau \in K\}.\]

One interpretation of these definitions is via ordered sequences.
Let $K$ be a simplicial complex on the index set $[n]$. The
vertices are ordered by their integer labels. If
$\sigma$ is a simplex of $K$ on vertices $\{i_{1},\ldots,i_{m}\}$
where $1\leq i_{1}<\cdots <i_{m}\leq n$, identify $\sigma$ with
the sequence $(i_{1},\ldots,i_{m})$. The simplices of $K$ are then
ordered left lexicographically.

Now suppose $K$ is a shifted complex. If $(i_{1},\ldots,i_{m})$ is
an $(m-1)$-dimensional simplex of $K$, then $K$ must contain every
simplex of dimension $m-1$ which is lexicographically less than
$(i_{1},\ldots,i_{m})$. Let $\mbox{rest}\{2,\ldots,n\}$ be the
simplicial subcomplex of $K$ which is defined as the collection of
simplices $(i_{1},\ldots,i_{m})\in K$ with $i_{1}\geq 2$. Observe
that $\mbox{star}(1)$ consists of those simplices
$(i_{1},\ldots,i_{m})$ for which
$(1,\ldots,i_{1}-1,i_{1},\ldots,i_{m})$ is also a simplex of $K$,
and $\mbox{link}(1)$ consists of those simplices which are in both
$\mbox{star}(1)$ and $\mbox{rest}\{2,\ldots,n\}$.

All this can now be formulated topologically in terms of coordinate
subspaces. Assume that $K$ is a simplicial complex on the index set $[n]$.
Let $X_{1},\ldots,X_{n}$ be path-connected spaces. Then we can associate
a coordinate subspace
$A$ of $X_{1}\times\cdots\times X_{n}$ to $K$ by letting $A$ be the union
of all subspaces $X_{i_{1}}\times\cdots\times X_{i_{m}}$ where
$(i_{1},\ldots,i_{m})$ is a simplex of $K$. Now suppose $K$ is shifted.
Let $\mbox{Star}(1)$, $\mbox{Link}(1)$, and $\mbox{Rest}\{2,\ldots,n\}$
be the coordinate subspaces of $X_{1}\times\cdots\times X_{n}$ associated
to $\mbox{star}(1)$, $\mbox{link}(1)$, and $\mbox{rest}\{2,\ldots,n\}$
respectively.

We now give a set of conditions on the inclusion of one shifted complex
into another which guarantees the existence of a regular sequence
between their corresponding coordinate subspaces.

\begin{proposition}
   \label{regseqs}
   Let $L$ and $K$ be two shifted complexes on the index set $[n]$,
   where $L$ is contained within $\mbox{star}(1)$ of
   $K$ and $K$ has no disjoint points. Fix path-connected spaces
   $X_{1},\ldots,X_{n}$. Let $B$ and $A$
   be the coordinate subspaces of $X_{1}\times\cdots\times X_{n}$ which
   correspond to $L$ and $K$ respectively. Let $\mbox{Star}(1)\subseteq A$
   be the coordinate subspace which corresponds to $\mbox{star}(1)\subseteq K$.
   Then there is a sequence of coordinate subspaces
   \[B=A_{0}\subseteq A_{1}\subseteq\cdots\subseteq A_{l}=\mbox{Star}(1)\]
   which is regular.
\end{proposition}

Before beginning with the proof of Proposition~\ref{regseqs} we give two
two examples which will be used subsequently. Observe that since all $n$
vertices are in $L$, the coordinate subspace $X_{1}\vee\cdots\vee X_{n}$
is contained in $B$.

\begin{example}
   \label{regex1}
Let $K$ be a connected shifted complex. Let $L$ be the disjoint
union of the $n$ vertices of $K$. Consider $\mbox{star}(1)$ in $K$.
Then $B=X_{1}\vee\cdots\vee X_{n}$, $A=\mbox{Star}(1)$, and
Proposition~\ref{regseqs} says that there exists a sequence of
coordinate subspaces
\[X_{1}\vee\cdots\vee X_{n}=A_{0}\subseteq A_{1}\subseteq\cdots
    \subseteq A_{l}=\mbox{Star}(1)\]
which is regular.
\end{example}

\begin{example}
   \label{regex2}
Let $K$ be a connected shifted complex. Let $L$ be $\mbox{link}(1)$.
Let $\mbox{star}_{R}(2)$ be $\mbox{star}(2)$ in
$\mbox{rest}\{2,\ldots,n\}$. Then $B=\mbox{Link}(1)$ and
$A=\mbox{Star}_{R}(2)$. To apply Proposition~\ref{regseqs} we need
to check that (within $\mbox{Rest}\{2,\ldots,n\}$) $\mbox{link}(1)$
is contained in $\mbox{star}_{R}(2)$. Let $(i_{1},\ldots,i_{m})$ be
a simplex of $\mbox{link}(1)$. If $i_{1}=2$ then
$(i_{1},\ldots,i_{m})$ is clearly in $\mbox{star}_{R}(2)$. If
$i_{1}\neq 2$, then as $\mbox{link}(1)\subseteq\mbox{star}(1)$ (in
$K$), the definition of $\mbox{link}(1)$ says there exists a simplex
$(1,2,\ldots,i_{1}-1,i_{1},\ldots,i_{m})$ in $\mbox{star}(1)$. The
restricted simplex $(2,\ldots,i_{1}-1,i_{1},\ldots,i_{m})$ is
therefore in $\mbox{star}_{R}(2)$. Thus $\mbox{star}_{R}(2)$
contains all the simplices in $\mbox{link}(1)$.
Proposition~\ref{regseqs} then implies that there is a sequence of
coordinate subspaces
\[\mbox{Link}(1)=A_{0}\subseteq A_{1}\subseteq\cdots
    \subseteq A_{l}=\mbox{Star}_{R}(1)\]
which is regular.
\end{example}

\noindent
\textit{Proof of Proposition~\ref{regseqs}}:
We adjoin subspaces to $B$ in two separate iterations. These adjunctions
correspond to gluing simplices to $L$ one at a time until $\mbox{star}(1)$
in $K$ is obtained.

\noindent \textbf{Iteration 1}: Since $K$ is connected and shifted,
every vertex in $K$ is connected by an edge to the vertex $1$, that
is, the simplex $(1,j)$ is in $K$ for every $2\leq j\leq n$. Now $L$
may contain disjoint points. If so, since $L$ is shifted, the
simplices $(1,j)$ will not be in $L$ for $j\geq j_{0}$, where
$j_{0}$ is the first vertex not connected to $1$. In terms of
coordinate subspaces, each $X_{j}$ is a wedge summand of $B$, and
$B$ contains the coordinate subspaces $X_{1}\times X_{j}$ for
$j<j_{0}$. The point of this first iteration is to adjoin the
coordinate subspaces $X_{1}\times X_{j}$ for $j\geq j_{0}$. They
will be adjoined in left lexicographical order. The adjunction is
realised by a homotopy pushout
\[\diagram
      X_{1}\vee X_{j}\rto\dto & A_{k-1}\dto \\
      X_{1}\times X_{j}\rto & A_{k}
  \enddiagram\]
which defines the space $A_{k}$. Here, we begin with the $j_{0}$
case, where $A_{0}=B$, so $k=j-1-j_{0}$. To show that this sequence
is regular, we need to show that there is a homotopy pushout
\[\diagram
      M_{k-1}\rto\dto & N_{k-1}\dto \\
      X_{1}\vee X_{j}\rto & A_{k-1}.
  \enddiagram\]
Take $M_{k-1}=X_{1}$. Observe that by the iteration to this point,
$A_{k-1}$ is the wedge $X_{1}\vee\cdots\vee X_{n}$ with the coordinate
subspaces $X_{1}\times X_{i}$ adjoined for $2\leq i\leq j-1$. In
particular, $X_{j}$ is a wedge summand of $A_{k-1}$. Let $N_{k-1}$
be the complementary wedge summand of $A_{k-1}$, so
$A_{k-1}\simeq X_{j}\vee N_{k-1}$. Then it is clear that
$M_{k-1}=X_{1}$ includes into $N_{k-1}$, the diagram
above homotopy commutes, and it is in fact a homotopy pushout.

\noindent
\textbf{Iteration 2}:
First observe that at the end of Iteration~$1$, all the coordinate
subspaces $X_{1}\times X_{j}$ for $2\leq j\leq n$
have  been adjoined to $\mbox{Star}(1)$. So
$A_{n-j_{0}}=X_{1}\times (X_{2}\vee\cdots\vee X_{n})$.

We now adjoin the remaining coordinate subspaces of $\mbox{Star}(1)$ in a
two-step process. The idea is to adjoin all the remaining coordinate
subspaces corresponding to the two-dimensional simplices
of $\mbox{star}(1)$ in lexicographic order, then the coordinate subspaces
corresponding to the three-dimensional simplices of $\mbox{star}(1)$,
and so on. Suppose all the coordinate subspaces corresponding to the
$(m-2)$-dimensional simplices in $\mbox{star}(1)$ have been adjoined.
Suppose $(1,i_{2},\ldots,i_{m})$ is the simplex of dimension $m-1$ of
least lexicographic order whose corresponding coordinate subspace
has not already been adjoined.
To perform the adjunction it is necessary that the coordinate subpspaces
corresponding to the boundary of
$(1,i_{2},\ldots,i_{m})$ have already been adjoined. The boundary is
composed of the simplices
\[(1,i_{2},\ldots,i_{m-1}),
    (1,i_{3},\ldots,i_{m}),\ldots, (1,i_{2},\ldots,i_{m-2},i_{m}),
    \ \mbox{and}\ (i_{2},\ldots,i_{m}).\]
All the coordinate subspaces corresponding to boundary simplices starting
with the vertex~$1$ have already been adjoined by inductive hypothesis: all
the simplices are of dimension~$m-2$
and are all clearly in $\mbox{star}(1)$. The lexicographical ordering
implies that the coordinate subspace corresponding to the simplex
$(i_{2},\ldots,i_{m})$ has not yet been adjoined. So we first need to
adjoin the coordinate subspace corresponding to $(i_{2},\ldots,i_{m})$
and then adjoin the coordinate subspace corresponding to
$(1,i_{2},\ldots,i_{m})$.
Note that the coordinate subspaces corresponding to the boundary simplices
of $(i_{2},\ldots,i_{m})$ have already been
adjoined because $\mbox{star}(1)$ being shifted means that if $\tau$
is a simplex in the boundary of $(i_{2},\ldots,i_{m})$ then $(1,\tau)$
is also a simplex of $\mbox{star}(1)$, and as its dimension is $m-2$,
the corresponding coordinate subspace has already been adjoined by
inductive hypothesis.

The two-step gluing process is realised by the homotopy pushouts
\[\diagram
      FW(i_{2},\ldots,i_{m})\rto\dto & A_{k-1}\dto
          & & FW(1,i_{2},\ldots,i_{m})\rto\dto & A_{k}\dto \\
      X_{i_{1}}\times\cdots\times X_{i_{m}}\rto & A_{k}
          & & X_{1}\times X_{i_{2}}\times\cdots\times X_{i_{m}}\rto & A_{k+1}
  \enddiagram\]
where the pushouts define the spaces $A_{k}$ and $A_{k+1}$. Observe that
if we assume $A_{k-1}\simeq X_{1}\times A^{\prime}_{k-1}$ -- this is true
for the base case $A_{n-j_{0}}$ as mentioned at the beginning of this
iteration -- then the two-step process in adjoining the coordinate
subspace corresponding to the simplex $(1,i_{2},\ldots,i_{m})$ implies
that $A_{k+1}\simeq X_{1}\times A^{\prime}_{k+1}$. Thus if we show that
the two-step process is itself a regular sequence, then the entire
iteration is a string of two-step regular sequences and so is a regular
sequence, completing the proof.

For the $k-1$ case, as $A_{k-1}\simeq X_{1}\times A^{\prime}_{k-1}$,
the definition of a regular sequence enforces us to project onto
$A^{\prime}_{k-1}$ and look for a homotopy pushout
\[\diagram
      M_{k-1}\rto\dto & N_{k-1}\dto \\
      FW(i_{2},\ldots,i_{m})\rto & A^{\prime}_{k-1}
  \enddiagram\]
where $M_{k-1}$ is a proper coordinate subspace of
$FW(i_{2},\ldots,i_{m})$. Having projected away from variable~$1$,
this homotopy pushout is really a lower dimensional case which
builds up $\mbox{Star}(2)$ within $\mbox{Rest}\{2,\ldots,n\}$. The
inductive hypothesis on dimension means that we can assume that this
homotopy pushout exists. For the $k$ case, we need to show that
there is a homotopy pushout
\[\diagram
      M_{k}\rto\dto & N_{k}\dto \\
      FW(1,i_{2},\ldots,i_{m})\rto & A{k}
  \enddiagram\]
where $M_{k}$ is a proper coordinate subspace of $FW(1,i_{2},\ldots,i_{m})$.
Let $M_{k-1}=X_{i_{2}}\times FW(i_{3},\ldots,i_{m})$. (Note that
$M_{k-1}$ equals $\mbox{Star}(1)$ in $FW(i_{2},\ldots,i_{m})$.) Observe
that if such a homotopy pushout exists, then
$N_{k}$ needs to contain all the coordinate subspaces of
$A_{k}$ except $X_{i_{2}}\times\cdots\times X_{i_{m}}$. But this is
exactly the description of $A_{k-1}$, so by taking $N_{k}=A_{k-1}$
we obtain the desired homotopy pushout.
$\qqed$

\section{The homotopy type of $\mathcal{Z}_{K}$ for shifted complexes}
\label{sec:ZKtype}

Recall that if $K$ is a simplicial complex on the index set $[n]$,
then there is a corresponding Davis-Januszkiewicz space $DJ(K)$ and
a homotopy fibration
\[\nameddright{\mathcal{Z}_{K}}{}{DJ(K)}{}{\prod_{i=1}^{n} BT}.\]
One of the main goals of the paper is to prove Theorem~\ref{shifted_wedge},
which we restate as:

\begin{theorem}
   \label{shifted}
   If $K$ is a shifted complex, then $\mathcal{Z}_{K}$ is homotopy
   equivalent to a wedge of spheres. That is, $K\in\mathcal{F}_{0}$.
\end{theorem}

It is well known (and easy to prove) that if $K$ is shifted then
each of $\mbox{link}(1)$, $\mbox{star}(1)$, and $\mbox{rest}\{2,\ldots,n\}$
is shifted, $\mbox{star}(1)=(1)\ast\mbox{link}(1)$, and there is a
topological pushout
\[\diagram
      \mbox{link}(1)\rto\dto & \mbox{rest}\{2,\ldots,n\}\dto \\
      \mbox{star}(1)\rto & K.
  \enddiagram\]
This results in a corresponding homotopy pushout of
Davis-Januszkiewicz spaces
\[\diagram
      DJ(\mbox{link}(1))\rto\dto & DJ(\mbox{rest}\{2,\ldots,n\})\dto \\
      DJ(\mbox{star}(1))\rto & DJ(K)
  \enddiagram\]
where $DJ(\mbox{star}(1))=BT\times DJ(\mbox{link}(1))$. Mapping the four
corners into $\prod_{i=1}^{n} BT$ and taking homotopy fibres gives a
cube as in Lemma~\ref{cube}, and in particular a homotopy pushout of fibres
\[\diagram
      \mathcal{Z}_{\mbox{link}(1)}\rto\dto
         & \mathcal{Z}_{\mbox{rest}\{2,\ldots,n\}}\dto \\
      \mathcal{Z}_{\mbox{star}(1)}\rto & \mathcal{Z}_{K}.
  \enddiagram\]
We wish to show that each of $\mathcal{Z}_{\mbox{link}(1)}$,
$\mathcal{Z}_{\mbox{rest}\{2,\ldots,n\}}$, and
$\mathcal{Z}_{\mbox{star}(1)}$ is homotopy equivalent to a wedge of
spheres, and then identify the maps in the homotopy pushout in order
to show that $\mathcal{Z}_{K}$ is also homotopy equivalent to a wedge
of spheres.

This topological problem can be reformulated more generally for coordinate
subspaces. We still assume that $K$ is a shifted complex on the
index set $[n]$. Let $X_{1},\ldots,X_{n}$ be path-connected spaces. Let
$A$ be the coordinate subspace of $X_{1}\times\cdots\times X_{n}$
associated to $K$. Then there is a homotopy pushout
\begin{equation}
   \label{coordApo}
   \diagram
       \mbox{Link}(1)\rto\dto & \mbox{Rest}\{2,\ldots,n\}\dto \\
       \mbox{Star}(1)\rto & A
   \enddiagram
\end{equation}
where $\mbox{Star}(1)\simeq X_{1}\times\mbox{Link}(1)$. Now compose each of
the four corners with the inclusion
\(\namedright{A}{}{X_{1}\times\cdots\times X_{n}}\)
and take homotopy fibres. Let $F_{L}$, $F_{S}$, $F_{R}$, and $F_{A}$ be
the homotopy fibres of the respective inclusions of $\mbox{Link}(1)$,
$\mbox{Star}(1)$, $\mbox{Rest}\{2,\ldots,n\}$, and $A$ into
$X_{1}\times\cdots\times X_{n}$. Then Lemma~\ref{cube} says there is
a homotopy pushout of fibres
\begin{equation}
   \label{FLpo}
   \diagram
        F_{L}\rto\dto & F_{R}\dto \\
        F_{S}\rto & F_{A}.
   \enddiagram
\end{equation}

The homotopy pushout in~\eqref{FLpo} can be refined. First, consider
the map \(\namedright{F_{L}}{}{F_{S}}\). As $\mbox{link}(1)$ is a
simplicial complex on the vertices $\{2,\ldots,n\}$, the space
$\mbox{Link}(1)$ is a coordinate subspace of
$X_{2}\times\cdots\times X_{n}$. Thus $F_{L}\simeq\Omega
X_{1}\times\overline{F}_{L}$ where $\overline{F}_{L}$ is the
homotopy fibre of the inclusion
\(\namedright{F_{L}}{}{X_{2}\times\cdots\times X_{n}}\). Continuing,
as $\mbox{Star}(1)\simeq X_{1}\times\mbox{Link}(1)$, there is a
homotopy pullback
\[\diagram
      \Omega X_{1}\times\overline{F}_{L}\rto\dto & \mbox{Link}(1)\rto\dto
         & X_{1}\times\cdots\times X_{n}\ddouble \\
      F_{S}\rto & X_{1}\times\mbox{Link}(1)\rto
         & X_{1}\times\cdots\times X_{n}.
  \enddiagram\]
As the map
\(\namedright{\mbox{Link}(1)}{}{X_{1}\times\mbox{Link}(1)}\)
is the inclusion of the second factor, the previous homotopy pullback
shows that $\overline{F}_{L}\simeq F_{S}$ and the map
\(\namedright{\Omega X_{1}\times\overline{F}_{L}}{}{F_{S}}\)
is the projection. Next, consider the map
\(\namedright{F_{L}}{}{F_{R}}\).
As $\mbox{Rest}\{2,\ldots,n\}$
is a coordinate subspace of $X_{2}\times\cdots\times X_{n}$, we have
$F_{R}\simeq\Omega X_{1}\times\overline{F}_{R}$ where $\overline{F}_{R}$
is the homotopy fibre of
\(\namedright{\mbox{Rest}\{2,\ldots,n\}}{}{X_{2}\times\cdots\times X_{n}}\).
As $\mbox{Link}(1)$ is a subspace of $\mbox{Rest}\{2,\ldots,n\}$, the map
\(\namedright{F_{L}}{}{F_{R}}\)
becomes
\(\llnamedright{\Omega X_{1}\times F_{S}}{\Omega X_{1}\times\gamma}
     {\Omega X_{1}\times\overline{F}_{R}}\)
for some map $\gamma$. Collecting all this information on the
homotopy fibres, the homotopy pushout in diagram~\eqref{FLpo}
becomes a homotopy pushout
\begin{equation}
   \label{refinedFLpo}
   \diagram
        \Omega X_{1}\times F_{S}\rto^-{\Omega X_{1}\times\gamma}\dto^{\pi_{2}}
           & \Omega X_{1}\times\overline{F}_{R}\dto \\
        F_{S}\rto & F_{A}.
   \enddiagram
\end{equation}
The goal is to identify the homotopy type of $F_{A}$. We do this
in Proposition~\ref{FAtype}. It may be useful to recall the definition
of $\mathcal{G}_{1}^{n}$ in~\ref{Gdefn}.

\begin{proposition}
   \label{FAtype}
   Let $K$ be a shifted complex on the index set $[n]$. Let
   $X_{1},\ldots,X_{n}$ be path-connected spaces and let $A$ be the
   coordinate subspace
   of $X_{1}\times\cdots\times X_{n}$ which corresponds to $K$.
   Use the notation and setup established in diagrams~\eqref{coordApo}
   and~\eqref{refinedFLpo}. Then the following hold:
   \begin{letterlist}
      \item $F_{S}\in\mathcal{G}_{1}^{n}$ and
            $\overline{F}_{R}\in\mathcal{G}_{2}^{n}$;
      \item there is a homotopy decomposition
            $F_{S}\simeq F_{S}^{1}\vee F_{S}^{2}$
            such that $F_{S}^{1},F_{S}^{2}\in\mathcal{G}_{1}^{n}$,
            the restriction of $\gamma$ to $F_{S}^{1}$ is null homotopic,
            and the restriction of $\gamma$ to $F_{S}^{2}$ has a right
            homotopy inverse;
      \item $F_{A}\in\mathcal{G}_{1}^{n}$.
   \end{letterlist}
\end{proposition}

\begin{proof}
We induct on $n$, the number of vertices. When $n=1$, we have
$A=X_{1}$, $\mbox{Star}(1)=X_{1}$, $\mbox{Rest}\{2,\ldots,n\}=\ast$,
and $\mbox{Link}(1)=\ast$. Composing into (the product space)
$X_{1}$ and taking homotopy fibres, we immediately see that
$F_{S}\simeq\ast$, $\overline{F}_{R}\simeq\ast$, $\gamma$ is
homotopic to the map from the basepoint to itself so part~(b)
trivially holds, and $F_{A}\simeq\ast$.

Assume the Proposition holds for $n-1$ vertices.
First, applying Proposition~\ref{F0toFl}~(a) to the regular sequence from
$X_{1}\vee\cdots\vee X_{n}$ to $\mbox{Star}(1)$ in Example~\ref{regex1}
shows that $F_{S}\in\mathcal{G}_{1}^{n}$. Next, since
$\mbox{Rest}\{2,\ldots,n\}$ is a shifted complex on the vertices
$\{2,\ldots,n\}$, the inductive hypothesis implies that
$\overline{F}_{R}\in\mathcal{G}_{2}^{n}$. This proves part~(a).

Assume part~(b) for the moment. Lemma~\ref{shiftedpo} then applies
to show there is a homotopy equivalence
\[F_{A}\simeq F^{2}_{S}\vee (\Omega X_{1}\ast F^{1}_{S}).\]
Thus $A\in\mathcal{G}_{1}^{n}$, proving part~(c).

To prove part~(b), we need to closely examine the map
\(\namedright{F_{S}}{\gamma}{\overline{F}_{R}}\). This was defined
in the setup for diagram~\eqref{refinedFLpo} by a homotopy pullback
\[\diagram
      F_{S}\rto\dto^{\gamma} & \mbox{Link}(1)\rto\dto
         & X_{2}\times\cdots\times X_{n}\ddouble \\
      \overline{F}_{R}\rto & \mbox{Rest}\{2,\ldots,n\}\rto
         & X_{2}\times\cdots\times X_{n}.
  \enddiagram\]
By definition, $\mbox{Star}_{R}(2)$ is a coordinate subspace of
$\mbox{Rest}\{2,\ldots,n\}$. In Example~\ref{regex2} we showed that
$\mbox{Link}(1)$ is a coordinate subspace of $\mbox{Star}_{R}(2)$.
Thus there is a diagram of iterated homotopy pullbacks
\[\diagram
      F_{S}\rto\dto^{\delta} & \mbox{Link}(1)\rto\dto
         & X_{2}\times\cdots\times X_{n}\ddouble \\
      F_{\overline{S}}\rto\dto^{\epsilon} & \mbox{Star}_{R}(2)\rto\dto
         & X_{2}\times\cdots\times X_{n}\ddouble \\
      \overline{F}_{R}\rto & \mbox{Rest}\{2,\ldots,n\}\rto
         & X_{2}\times\cdots\times X_{n}
  \enddiagram\]
where the pullbacks define the space $F_{\overline{S}}$ and
the maps $\delta$ and $\epsilon$. Hence $\gamma\simeq\epsilon\circ\delta$.
We deal with each of $\delta$ and $\epsilon$ one at a time.

Applying Proposition~\ref{F0toFl}~(b) to the regular sequence from
$\mbox{Link}(1)$ to $\mbox{Star}_{R}(2)$ in Example~\ref{regex2} shows that
$F_{S}\simeq E_{1}\vee E_{2}$ where $E_{1},E_{2}\in\mathcal{G}_{1}^{n}$,
the restriction of $\delta$ to $E_{1}$ is null homotopic, and the
restriction of $\delta$ to $E_{2}$ has a right homotopy inverse.

For $\epsilon$, we appeal to the inductive hypothesis. Let
$\mbox{link}_{R}(2)$ be $\mbox{link}(2)$ within $\mbox{rest}\{2,\ldots,n\}$.
Since $\mbox{rest}\{2,\ldots,n\}$ is a shifted complex, it is the pushout
of $\mbox{star}_{R}(2)$ and $\mbox{rest}\{3,\ldots,n\}$ over
$\mbox{link}_{R}(2)$. This results in a homotopy pushout of the
corresponding coordinate subspaces (in $X_{2}\times\cdots\times X_{n}$)
\[\diagram
     \mbox{Link}_{R}(2)\rto\dto & \mbox{Rest}\{3,\ldots,n\}\dto \\
     \mbox{Star}_{R}(2)\rto & \mbox{Rest}\{2,\ldots,n\}.
  \enddiagram\]
Let $F_{\overline{L}}$ and $F_{\overline{R}}$ respectively be the
homotopy fibres of the inclusions of $\mbox{Link}_{R}(2)$ and
$\mbox{Rest}\{3,\ldots,n\}$ into $X_{2}\times\cdots\times X_{n}$.
Recall that $F_{\overline{S}}$ and $\overline{F}_{R}$ respectively
have been defined as the homotopy fibres of the inclusions of
$\mbox{Star}_{R}(2)$ and $\mbox{Rest}\{2,\ldots,n\}$ into
$X_{2}\times\cdots\times X_{n}$. As in diagram~\eqref{FLpo}, when
all four corners of the pushout above are mapped into
$X_{2}\times\cdots\times X_{n}$, we obtain a homotopy pushout of
fibres
\[\diagram
     F_{\overline{L}}\rto\dto & F_{\overline{R}}\dto \\
     F_{\overline{S}}\rto & \overline{F}_{R}.
  \enddiagram\]
Arguing as for diagram~\eqref{refinedFLpo}, this homotopy pushout of
fibres refines to a homotopy pushout
\[\diagram
     \Omega X_{2}\times F_{\overline{S}}\rto^-{1\times\overline{\gamma}}\dto
        & \Omega X_{2}\times\overline{F}_{\overline{R}}\dto \\
     F_{\overline{S}}\rto & \overline{F}_{R}
  \enddiagram\]
where $\overline{F}_{\overline{R}}$ is the homotopy fibre of the
inclusion
\(\namedright{\mbox{Rest}\{3,\ldots,n\}}{}{X_{3}\times\cdots\times
X_{n}}\). Let
\(\varphi:\namedright{F_{\overline{S}}}{}{\overline{F}_{R}}\) be the
map along the bottom row. Since the underlying shifted complex
$\mbox{rest}\{2,\ldots,n\}$ is on $n-1$ vertices, by inductive
hypothesis we can assume that there is a homotopy decomposition
$F_{\overline{S}}\simeq D_{1}\vee D_{2}$ where
$D_{1},D_{2}\in\mathcal{G}_{2}^{n}$, the restriction of
$\overline{\gamma}$ to $D_{1}$ is null homotopic while the
restriction of $\overline{\gamma}$ to $D_{2}$ has a right homotopy
inverse. Lemma~\ref{shiftedpo} then implies that there is a homotopy
equivalence $\overline{F}_{R}\simeq D_{2}\vee (\Omega X_{2}\ast
D_{1})$, the restriction of $\varphi$ to $D_{1}$ is null homotopic
while the restriction of $\varphi$ to $D_{2}$ has a right homotopy
inverse.

Now consider the composite
\(\gamma:\nameddright{F_{S}}{\delta}{F_{\overline{S}}}{\epsilon}
     {\overline{F}_{R}}\).
Since the restriction of $\delta$ to $E_{2}$ has a right homotopy
inverse, the decomposition $F_{\overline{S}}\simeq D_{1}\vee D_{2}$
results in a decomposition $E_{2}\simeq E^{1}_{2}\vee E^{2}_{2}$
where $E^{i}_{2}$ retracts off $D_{i}$. Let $F_{S}^{1}=E^{1}\vee E_{2}^{1}$,
and let $F_{S}^{2}=E_{2}^{2}$. Then the conditions on $\delta$
and $\epsilon$ imply that the restriction of $\gamma$ to $F_{S}^{1}$
is null homotopic while the restriction of $\gamma$ to $F_{S}^{2}$
has a right homotopy inverse, proving part~(b).
\end{proof}

With Proposition~\ref{FAtype} in hand, we can prove
Theorem~\ref{shifted} as a special case.
\medskip

\noindent
\textit{Proof of Theorem~\ref{shifted}}: In this case, each space
$X_{i}$ equals $BT$, the classifying space of the torus, the coordinate
subspace $A$ equals $DJ(K)$, and the homotopy fibre $F_{A}$ equals
$\mathcal{Z}_{K}$. Proposition~\ref{FAtype}~(c) says that
$\mathcal{Z}_{K}\in\mathcal{G}_{1}^{n}$, meaning that
$\mathcal{Z}_{K}$ is homotopy equivalent to a wedge of summands of
the form $\Omega BT_{i_{1}}\ast\cdots\ast\Omega BT_{i_{m}}$. Such a
summand is homotopy equivalent to $S^{m+1}$ since
$\Omega BT\simeq S^{1}$. Thus $\mathcal{Z}_{K}$ is homotopy equivalent
to a wedge of spheres, and so $K\in\mathcal{F}_{0}$.
$\qqed$
\medskip

A special case of a shifted complex is the full $i$-skeleton
$\Delta^{i}(n)$ of the standard simplex $\Delta^{n-1}$ on $n$
vertices. For path-connected spaces $X_{1},\ldots,X_{n}$, let
$T^{n}_{k}$ be the coordinate subspace associated to
$\Delta^{n-k}(n)$. Specifically,
\[T^{n}_{k}=\{(x_{1},\ldots,x_{n})\in X_{1}\times\cdots\times X_{n}
    \, \vert\ \mbox{at least $k$ of the $x_{i}$'s are $\ast$}\}.\]
In particular, $T^{n}_{0}=X_{1}\times\cdots\times X_{n}$,
$T^{n}_{1}$ is the fat wedge, $T^{n}_{n-1}=X_{1}\vee\cdots\vee
X_{n}$, and $T^{n}_{n}=\ast$. Let $F^{n}_{k}$ be the homotopy fibre
of the inclusion \(\namedright{T^{n}_{k}}{}{T^{n}_{0}}\). By
Theorem~\ref{shifted}, $F^{n}_{k}$ is homotopy equivalent to a wedge
of spheres. This wedge can be calculated explicitly using the
iteration in Proposition~\ref{regseqs} to reproduce a result first
obtained in a different context by Porter~\cite{P}. For a space $X$
and a positive integer $j$, let $j\cdot X$ be the wedge sum of $j$
copies of $X$. Let $X^{(j)}$ be the $j$-fold smash of $X$ with
itself.

\begin{theorem}[Porter]
   \label{Porter}
   For $n\geq 1$, let $X_{1},\ldots, X_{n}$ be path-connected
   spaces. Let $k$ be such that $1\leq k\leq n-1$. Then there
   is a homotopy equivalence
   \[F^{n}_{k}\simeq\bigvee_{j=n-k+1}^{n}
      \left(\bigvee_{1\leq i_{1}<\cdots<i_{j}\leq n}\ \binom{j-1}{n-k}
      \Sigma^{n-k}\Omega X_{i_{1}}\wedge\cdots\wedge\Omega X_{i_{j}}\right).\]
   $\qqed$
\end{theorem}

\begin{corollary}
   \label{Portercor}
   As in Theorem~\ref{Porter}, if $X_{i}=X$ for each $1\leq i\leq n$
   then there is a homotopy equivalence
   \[F^{n}_{k}\simeq\bigvee_{j=n-k+1}^{n}\binom{n}{j}\binom{j-1}{n-k}
      \Sigma^{n-k} (\Omega X)^{(j)}.\]
   $\qqed$
\end{corollary}

The case of relevance to us is Corollary~\ref{Portercor} applied
when $X=\mathbb{C}P^{\infty}$. Then $T^{n}_{k}$ corresponds to
$DJ(K)$ for $K=\Delta^{n-k}(n)$. The homotopy fibre $F^{n}_{k}$
corresponds to $\mathcal{Z}_{K}$. Since $\Omega X\simeq S^{1}$, we
obtain:

\begin{corollary}
   \label{ZkforDelta}
   If $K=(\Delta^{n-k}(n)$, then
   \[\mathcal{Z}_{K}\simeq\bigvee_{j=n-k+1}^{n}
       \binom{n}{j}\binom{j-1}{n-k} S^{n-k+j}.\]
   $\qqed$
\end{corollary}

\section{Topological extensions}
\label{sec:topext}

At this point, we have shown that if a simplicial complex $K$ is
shifted then its moment-angle complex $\mathcal{Z}_{K}$ is homotopy
equivalent to a wedge of spheres. Next, we want to consider other
simplicial complexes $K$ for which $\mathcal{Z}_{K}$ is homotopy
equivalent to a wedge of spheres, or for which $\mathcal{Z}_{K}$ is
stably equivalent to a wedge of spheres. We phrase this as follows.

\begin{definition}
Let $\F_t$ be the family of simplicial complexes $K$ for which
the moment-angle complex $\Z_K$ has the property that $\Sigma^t
\Z_K$ is homotopy equivalent to a wedge of spheres.
\end{definition}

The definition of $\F_t$ gives rise to a filtration
\[
\F_0\subset
\F_1\subset\ldots\subset\F_t\subset\ldots\subset\F_\infty.
\]
This filtration does not account for all simplicial
complexes $K$. As mentioned in the introduction, torsion can occur in the
cohomology ring of $\Z_K$ for certain simplicial complexes $K$,
making it impossible for $\Z_K$ to be even stably homotopic to a
wedge of spheres.

The problem we want to consider next is how the filtration level is
affected when combinatorial operations are applied to two simplicial
complexes from possibly different filtration levels. The
combinatorial operations we look at are: the disjoint union of
simplicial complexes, gluing along a common face and the join of
simplicial complexes. Recall that for given simplicial complexes
$K_1$ and $K_2$ on sets $\S_1$ and $\S_2$ respectively the
\textit{join} $K_1*K_2$ is the simplicial complex
\begin{equation*}
K_1*K_2:=\{\sigma\subset\S_1\cup\S_2 \,\vert\
\sigma=\sigma_1\cup\sigma_2, \sigma_1\in K_1, \sigma_2\in K_1\}
\end{equation*}
on the set $\S_1\cup\S_2$.

\begin{theorem}
\label{families2} Let $K_{1}\in \F_{t}$ and $K_{2}\in %
\ensuremath{\mathcal{F}}_{s}$ for some non-negative integers $t$
and~$s$. The effect on family membership of the simplicial
complex $K$ resulting from the following operations on $K_1$ and $K_2$ is:

\begin{enumerate}
\item  gluing along a common face:\newline
if $K=K_{1}\bigcup_{\sigma }K_{2}$, then $K\in \F_{m}$ where
$\sigma$ is a common face of $K_{1}$ and $K_{2}$ and ${m=\max
\{t,s\}}$;
\item  the disjoint union of simplicial complexes:\newline
if $K=K_{1}\coprod K_{2}$, then $K\in \ensuremath{\mathcal{F}}_{m}$
where $m=\max \{t,s\} $;
\item the join of simplicial complexes:\newline
if $K=K_{1}\ast K_{2}$, then $K\in \F_{m}$ where $m=\max \{t,s\}+1$.
\end{enumerate}
\end{theorem}
\begin{proof}

\indent $\mathrm{(1)}$ Let $DJ(K_{i})$, $i=1,2$, $DJ(\sigma)$ and
$DJ(K)$ be the corresponding Davis-Januszkiewicz spaces. Each vertex
in $K_{i}$, $\sigma$ or $K$ corresponds to a coordinate in
$DJ(K_{i})$, $DJ(\sigma)$ or $DJ(K)$ respectively. List the vertices
of $K_{1}$ as $\{1,\ldots,l,\ldots,m\}$, where the vertices of
$\sigma$ are $\{l+1,\ldots,m\}$. List the vertices of $K_{2}$ as
$\{l+1,\ldots,m,\ldots,n\}$. Regard $DJ(K_{1})$ as a subspace of
$\prod_{i=1}^{m}\mathbb{C}P^{\infty}$. Let $D_{1}$ be the image of
$DJ(K_{1})$ under the map
\(\namedright{\prod_{i=1}^{m}\mathbb{C}P^{\infty}}{}
    {\prod_{i=1}^{n}\mathbb{C}P^{\infty}}\)
given by the inclusion of the first $m$ coordinates. Similarly,
regard $DJ(K_{2})$ as a subspace of
$\prod_{i=l+1}^{n}\mathbb{C}P^{\infty}$, and let $D_{2}$ be its
image under the map
\(\namedright{\prod_{i=l+1}^{n}\mathbb{C}P^{\infty}}{}
     {\prod_{i=1}^{n}\mathbb{C}P^{\infty}}\)
given by the inclusion of the last $n-l$ coordinates. Since $\sigma$
is a simplex, $DJ(\sigma)$ is a product of $m-l$ copies of
$\mathbb{C}P^{\infty}$. Let $D_{3}$ be the image of $DJ(\sigma)$ in
$\prod_{i=1}^{n}\mathbb{C}P^{\infty}$ under the map
\(\namedright{\prod_{i=l+1}^{m}\mathbb{C}P^{\infty}}{}
     {\prod_{i=1}^{n}\mathbb{C}P^{\infty}}\)
given by the inclusion of the middle $m-l$ coordinates. Let $D$ be
the topological pushout
\begin{equation}\label{pushputD}
\xymatrix{
D_{3} \ar[d]\ar[r] & D_{1}\ar[d] \\
     D_{2}\ar[r] & D.}
  \end{equation}
Then $D=DJ(K)$ and is a subspace of
$\prod_{i=1}^{n}\mathbb{C}P^{\infty}$.

For notational convenience, let
$BT^{n}=\prod_{i=1}^{n}\mathbb{C}P^{\infty}$. Map each of the four
corners of pushout~\eqref{pushputD} into $BT^{n}$ and take homotopy
fibres. This gives homotopy fibrations
\[\nameddright{F}{}{D}{}{BT^{n}}\]
\[\nameddright{F_{1}\times N}{}{D_{1}}{}{BT^{n}}\]
\[\nameddright{M\times F_{2}}{}{D_{2}}{}{BT^{n}}\]
\[\nameddright{M\times N}{}{D_{3}}{}{BT^{n}}\]
where the first homotopy fibration defines $F$, $F_{1}$ is the
homotopy fibre of
\(\namedright{D_{1}}{}{\prod_{i=1}^{m}\mathbb{C}P^{\infty}}\),
$F_{2}$ is the homotopy fibre of
\(\namedright{D_{2}}{}{\prod_{i=l+1}^{n}\mathbb{C}P^{\infty}}\),
$M=\prod_{i=1}^{l} S^{1}$, and $N=\prod_{i=m+1}^{n} S^{1}$.
Including $D_{3}$ into $D_{1}$ gives a homotopy pullback diagram
\[\diagram
     \Omega BT^{n}\rto\ddouble & M\times N\rto\dto^{\theta}
         & D_{3}\rto\dto & BT^{n}\ddouble \\
     \Omega BT^{n}\rto & F_{1}\times N\rto & D_{1}\rto & BT^{n}
  \enddiagram\]
for some map $\theta$ of fibres. We now identify $\theta$. With
$BT^{m}=\prod_{i=1}^{m}\mathbb{C}P^{\infty}$, the pullback just
described is the product of the homotopy pullback
\[\diagram
     \Omega BT^{m}\rto\ddouble & M\rto\dto^{\theta^{\prime}}
         & D_{3}\rto\dto & BT^{m}\ddouble \\
     \Omega BT^{m}\rto & F_{1}\rto & D_{1}\rto & BT^{m}
  \enddiagram\]
and the path-loop fibration
\(\nameddright{N}{}{\ast}{}{\prod_{i=m+1}^{n}\mathbb{C}P^{\infty}}\).
So $\theta=\theta^{\prime}\times N$. Further, $M=\prod_{i=1}^{l}
S^{1}$ is a retract of $\Omega BT^{m}\simeq\prod_{i=1}^{m} S^{1}$
and \(\namedright{\Omega BT^{m}}{}{F_{1}}\) is null homotopic
since $\Omega BT^{m}$ is a retract of $\Omega D_{1}=\Omega
DJ(K_{1})$. Hence $\theta^{\prime}\simeq\ast$ and so
$\theta\simeq\ast\times N$. A similar argument for the inclusion of
$D_{3}$ into $D_{2}$ shows that the map of fibres
\(\namedright{M\times N}{}{M\times F_{2}}\) is homotopic to
$M\times\ast$.

Collecting all this information about homotopy fibres,
Lemma~\ref{cube} shows that there is a homotopy pushout
\[\diagram
      M\times N\rto^-{\ast\times N}\dto^{M\times\ast}
         & F_{1}\times N\dto \\
      M\times F_{2}\rto & F.
  \enddiagram\]
Lemma~\ref{gluingpo} then gives a homotopy decomposition
\[F\simeq (M\ast N)\vee (M\ltimes F_{2})\vee (F_{1}\rtimes N).\]
We want to show that $\Sigma^{m} F$ is homotopy equivalent to a
wedge of spheres, where $m=\max\{t,s\}$. If so, then as $D=DJ(K)$,
we have $F=\mathcal{Z}_{K}$ and hence $K\in\mathcal{F}_{\max\{t,s\}}$,
proving~(1).

To show $\Sigma^{m} F$ is homotopy equivalent to a wedge of spheres,
we show that each of $\Sigma^{m}(M\ast N)$, $\Sigma^{m} (M\ltimes F_{2})$,
and $\Sigma^{m} (F_{1}\rtimes N)$ are homotopy equivalent to wedges
of spheres. First, observe that the suspension of a product of spheres
is homotopy equivalent to a wedge of spheres. Since $M$ and $N$ are
products of copies of $S^{1}$, $M\ast N$ is therefore homotopy equivalent
to a wedge of spheres, and hence $\Sigma^{m} (M\ast N)$ is as well.
Second, if $m\geq 1$ then
$\Sigma^{m} (M\ltimes F_{2})\simeq
   \Sigma^{m} (M\wedge F_{2})\vee\Sigma^{m} F_{2}$.
By hypothesis, $\Sigma^{m} F_{2}$ is homotopy equivalent to a wedge
of spheres. As $m\geq 1$, $\Sigma^{m} (M\wedge F_{2})$ is the $(m-1)$-fold
suspension of $(\Sigma M)\wedge F_{2}$. But $\Sigma M$ is homotopy
equivalent to a wedge of spheres, so $(\Sigma M)\wedge F_{2}$ is
homotopy equivalent to a wedge of suspensions of $F_{2}$. Therefore
$\Sigma^{m} (M\wedge F_{2})$ is homotopy equivalent to a wedge of
$m$-fold suspensions of $F_{2}$, implying that it is homotopy
equivalent to a wedge of spheres. If $m=0$,
then $F_{2}$ is still homotopy equivalent to a wedge of spheres, so
we can write $F_{2}\simeq\Sigma F_{2}^{\prime}$, where $F_{2}^{\prime}$
is a wedge of spheres. We then have
$M\ltimes F_{2}\simeq M\ltimes (\Sigma F_{2}^{\prime})
    \simeq \Sigma (M\ltimes F_{2}^{\prime})$
and the decomposition into a wedge of spheres now follows as in the
$m\geq 1$ case. The decomposition of the summand $F_{1}\rtimes N$
into a wedge of spheres is exactly as for $M\ltimes F_{2}$.

$\mathrm{(2)}$ Let $K=K_{1}\coprod K_{2}$ be the disjoint union of
two simplicial complexes $K_{1}$ and $K_{2}$ on the index sets $[m]$
and $[n]$ respectively. Then their disjoint union $K=K_1\coprod
K_2$ is a simplicial complex on the index set $[m+n]$ obtained as the
result of gluing $K_{1}$ to $K_{2}$ along the empty face. Applying
part~$\mathrm{(1)}$ then shows that $K\in
\ensuremath{\mathcal{F}}_{m}$, where $m=\max \{t,s\} $. Moreover, the
homotopy type of $\Z_K$ is given by
\[
\ensuremath{\mathcal{Z}}_K\simeq \Big(\prod^{m}_{i=1} S^1*
\prod^{n}_{j=1} S^1\Big) \vee \Big(\ensuremath{\mathcal{Z}}_{K_1} \rtimes
\prod^{n}_{i=1}S^1\Big) \vee \Big(\prod^{m}_{i=1}
S^1\ltimes\ensuremath{\mathcal{Z}}_{K_2}\Big).
\]

$\mathrm{(3)}$ The Stanley-Reisner ring of the join $K_1*K_2$ of two
simplicial complexes $K_1$ and $K_2$ on the index sets $[m]$ and
$[n]$ has the following form:
\begin{equation*}
\mathbb{Z}[K_1*K_2]=\mathbb{Z}[K_1]\otimes \mathbb{Z}[K_2].
\end{equation*}
Therefore the fibration
\begin{equation*}
DJ(K_1*K_2)\longrightarrow BT^{m+n}
\end{equation*}
associated to the join of $K_1$ and $K_2$ is the product fibration
\begin{equation*}
DJ(K_1)\times DJ(K_2)\longrightarrow BT^{m}\times BT^{n}.
\end{equation*}
Hence $\ensuremath{\mathcal{Z}}_{K_1*K_2}\simeq\ensuremath{\mathcal{Z}}%
_{K_1}\times\ensuremath{\mathcal{Z}}_{K_2}$. This proves part
$\mathrm{(3)}$ and finishes the proof of the Theorem.
\end{proof}
As a corollary we specify the operations on simplicial complexes for
which $\F_{0}$ is closed.

\begin{corollary}
Let $K_{1}$ and $K_{2}$ be simplicial complexes in $\ensuremath{%
\mathcal{F}}_{0}$. Then $\F_{0}$ is closed for the following
operations on simplicial complexes:

\begin{enumerate}
\item  gluing along a common face,\newline
$K=K_{1}\bigcup_{\sigma }K_{2}\in \F_{0}$, where $%
\sigma $ is a common face of $K_{1}$ and $K_{2}$.
\item  the disjoint union of simplicial complexes,\newline
$K=K_{1}\coprod K_{2}\in \F_{0}$;
\end{enumerate}
\end{corollary}

\section{Algebra}
\label{sec:algebra}

Let $A$ be a polynomial ring on $n$ variables $k[x_1,\ldots,x_n]$
over a field $k$ and let $R=A/I$, where $I$ is homogeneous ideal. In
this section we shall be interested in the nature of $\Tor _R(k,k)$;
specifically, in identifying a class of rings $R$ for which all
Massey products in $\Tor_A(R,k)$ vanish and how this impacts upon
the Poincar\'{e} series of $R$. Recall that the Poincar\'{e} series
of $R$ is the formal power series
\[
P(R)=\sum^{\infty}_{i=0}b_it^i
\]
where $b_i=\dim_k \Tor_R^i(k,k)$ are the Betti numbers of $R$. It
has been conjectured by Kaplansky and Serre that $P(R)$ always
represents a rational function. The regular local rings were the
first rings for which $P(R)$ was explicitly computed. In this case
Serre~\cite{Se} showed that $P(R)=(1+t)^n$. Tate~\cite{T} showed
that if $R$ is a complete intersection, then there exist
non-negative integers $m,n$ such that
\[
P(R)=\frac{(1+t)^n}{(1-t^2)^m}.
\]
Golod~\cite{Go} made a far reaching contribution to the problem by
showing that if certain homology operations on the Koszul complex
vanish, then there exist non-negative integers $n,c_1,\ldots ,c_n$
such that
\[
P(R)=\frac{(1+t)^n}{1-\sum^n_{i=1}c_it^{i+1}}.
\]
In general not much is known about the rationality of $P(R)$;
although there is an inequality due to Golod~\cite{Go} showing that
$P(R)$ is always bounded (coefficient-wise) by a rational function.

In the past, describing various properties of $\Tor_R(k,k)$ has been
largely an algebraic problem. Further on, we translate the problem
of rationality of the Poincar\' e series into topology by using
recent results of toric topology. Then by using our results on the
homotopy type of the complement of a coordinate subspace arrangement,
we find a class of rings $R$ for which $P(R)$ is a rational function
determined by $P(\Tor_A(R,k))$ .

In what follows $R$ will be the Stanley-Reisner ring $k[K]$ of an
arbitrary simplicial complex $K$ on $n$ vertices. Recall from
Definition~\ref{Goloddef} that
the Stanley-Reisner ring $k[K]$ is {\it Golod} if all Massey
products in $\Tor_{k[v_1,\ldots,v_n]}(k[K],k)$ vanish.
Buchstaber and Panov~\cite{BP} proved that
\[
\Tor^*_{k[K]}(k,k)\cong H^*(\Omega DJ(K);k).
\]
This isomorphism now lets us exploit the topological properties
of the loop space $\Omega DJ(K)$ to obtain further information
about $\Tor_R(k,k)$. Looking at the split fibration
\[
\ddr{\Omega\Z_K}{}{\Omega DJ(K)}{}{T^n}
\]
we have
\[
\Tor_R^* (k,k)\cong H^*(\Omega DJ(K))=H^*(T^n)\otimes H^*(\Omega
\Z_K).
\]

A calculation using the bar resolution shows that
\[
P(H^*(\Omega \Z_K))\leq P(T(\Sigma^{-1} H^*(\Z_K)))
\]
where $\Sigma^{-1} H^*(\Z_K)$ is the desuspension of the module
$H^*(\Z_K)$. Therefore
\begin{equation*}
\label{golodin}
P(R)\leq (1+t)^n P(T(\Sigma^{-1} H^*(\Z_K))=
   \frac{t(1+t)^n}{t-P(H^*(\Z_K))}.
\end{equation*}

Looking at the Eilenberg-More spectral sequence (the bar resolution)
that computes the cohomology of the fibre in the path-loop fibration
$\ddr{\Omega\Z_K}{}{*}{}{\Z_K}$, we conclude that the above equality
is reached when the differentials are trivial. According to May, the
differentials are determined by the Massey products and therefore
they are trivial when all the Massey products in $H^*(\Z_K)$ vanish.
As $H^*(\Z_K)\cong \Tor_{k[v_1,\ldots,v_n]}(k[K],k)$~\cite{BP}, an
equality for $P(R)$ is obtained when the Stanley-Reisner ring $k[K]$
is Golod. This proves the following theorem.
\begin{theorem}
For a simplicial complex $K$,
\begin{equation}
\label{poincare} P(k[K])\leq\frac{t(1+t)^n}{t-P(H^*(\Z_K))}.
\end{equation}
Equality is obtained when $k[K]$ is Golod.
\end{theorem}
We proceed by describing a new class of Golod rings using
topological methods.
\begin{theorem}
\label{golodclass} If $K\in \mathcal{F}_{0}$, then $k[K]$ is a Golod
ring.
\end{theorem}
\begin{proof}
By definition of the family $\mathcal{F}_0$, when $K\in
\mathcal{F}_0$ then $\Z_K$ is homotopy equivalent to a wedge of
spheres. Therefore in the cohomology of $\Z_K$ all cup products and
higher Massey products are trivial. On the other hand, recall that
Buchstaber and Panov~\cite{BP} proved that
\[
H^*(\Z_K)\cong \Tor_{k[v_1,\ldots,v_n]}(k[K],k).
\]
Therefore in $\Tor_{k[v_1,\ldots,v_n]}(k[K],k)$ all Massey products
are trivial. Now by definition, the ring $k[K]$ is Golod.
\end{proof}
We finish by proving that the Poincar\' e series of a ring belonging
to the class defined in Theorem~\ref{golodclass} represents a
rational function.
\begin{corollary}
If $K\in \mathcal{F}_0$, then the Poincar\' e series of the ring
$k[K]$ has the following form
\[
P(k[K])=\frac{t(1+t)^n}{t-P(H^*(\Z_K))}.
\]
\end{corollary}
\begin{proof}
As $k[K]$ is a Golod ring, in \eqref{poincare} equality holds.
\end{proof}

\bibliographystyle{amsalpha}

\end{document}